\documentclass[a4paper,fleqn]{cas-sc}

\usepackage[square,numbers,sort&compress]{natbib}

\usepackage[figurename=Fig.]{caption}
\usepackage[title]{appendix}
\usepackage{subfig}
\usepackage{graphicx}
\usepackage{float}
\usepackage[capitalise]{cleveref}
\graphicspath{{figs/}}
\usepackage{mathrsfs} 
\usepackage[scr=esstix,cal=boondox]{mathalfa} 
\usepackage{tikz}
\usetikzlibrary{positioning, shapes.geometric}
\renewcommand{\printorcid}{} 
\ExplSyntaxOn              
\keys_set:nn { stm / mktitle } { nologo }   
\ExplSyntaxOff    

\def\tsc#1{\csdef{#1}{\textsc{\lowercase{#1}}\xspace}}
\tsc{WGM}
\tsc{QE}
\tsc{EP}
\tsc{PMS}
\tsc{BEC}
\tsc{DE}

\begin{document}
\let\WriteBookmarks\relax
\def\floatpagepagefraction{1}
\def\textpagefraction{.001}
\shorttitle{}
\shortauthors{Q. Liu et~al.}

\title [mode = title]{Adaptive coupling peridynamic least-square minimization with finite element method for fracture analysis}                      

\author[1]{Qibang Liu}[type=editor,
                        auid=000,bioid=1,
                        prefix=,
                        role=,
                        orcid=]
\cormark[1]
\ead{qibangliu@ksu.edu.}

\address[1]{Department of Mechanical and Nuclear Engineering, Kansas State University, Manhattan, KS 66506, USA}
\address[2]{Parks College of Engineering, Aviation, and Technology, Saint Louis University, MO 63103,USA}
\author[1]{X.J. Xin}
\author[2]{Jeff Ma}
\cortext[cor1]{Corresponding author. }

\begin{abstract}
This study presents an adaptive coupling peridynamic least-square minimization with the finite element method  (PDLSM-FEM) for fracture analysis. The presented method utilizes the PDLSM modeling discontinuities while maximizing the FEM region for computational efficiency. Within the presented adaptive PDLSM-FEM, only elements intersecting with the crack path and their neighboring elements are defined as PD elements, whose stiffness matrices are derived based on PDLSM equations. The remaining elements are conventional finite elements. Numerical integration of interaction integral is proposed and implemented to evaluate the stress intensity factors (SIFs) for 2-D problems. The criterion of maximum hoop tensile stress is employed for failure prediction. New contributions of this work include the adaptive coupling of PDLSM with FEM for minimizing the PD region and the application of the adaptive PDLSM-FEM to quasi-static crack propagation analysis. Simulations of three 2-D plane stress plates and one 3-D block with static or quasi-static cracks propagation are performed. Results show the proposed method improves computational efficiency substantially and has reasonable accuracy and good capability of crack propagation prediction.

\end{abstract}


\begin{keywords}
Peridynamics \sep FEM  \sep Weighted residual  \sep Fracture \sep Adaptive
\end{keywords}

\maketitle 

\section{Introduction}

The finite element method has been one of the most popular and successful numerical tools for studying structure behaviors for decades. 
The FEM requires the domain to be discretized into non-overlapping elements whose displacements are approximated with polynomials. 
Stress and strain fields then are derived from the displacement approximation based on classical continuum theory. 
Despite its widespread applications, the FEM suffers from drawbacks in handling displacement discontinuities, as in the cases of cracked bodies and other defects, since the FEM is formulated within the framework of classical continuum mechanics, in which the governing equations require the derivatives of displacements and become invalid at displacement discontinuities. Various ways to handle displacement discontinuities such as cracks have been proposed. 
A typical technique is to remesh repeatedly to match the discontinuities so that the discontinuities, such as cracks, always coincide with element edges and propagate only along the boundaries of neighboring elements. Remeshing, however, is computationally intensive and somewhat cumbersome.

 \begin{figure}[htbp]
    \centering
    \includegraphics[width=4in]{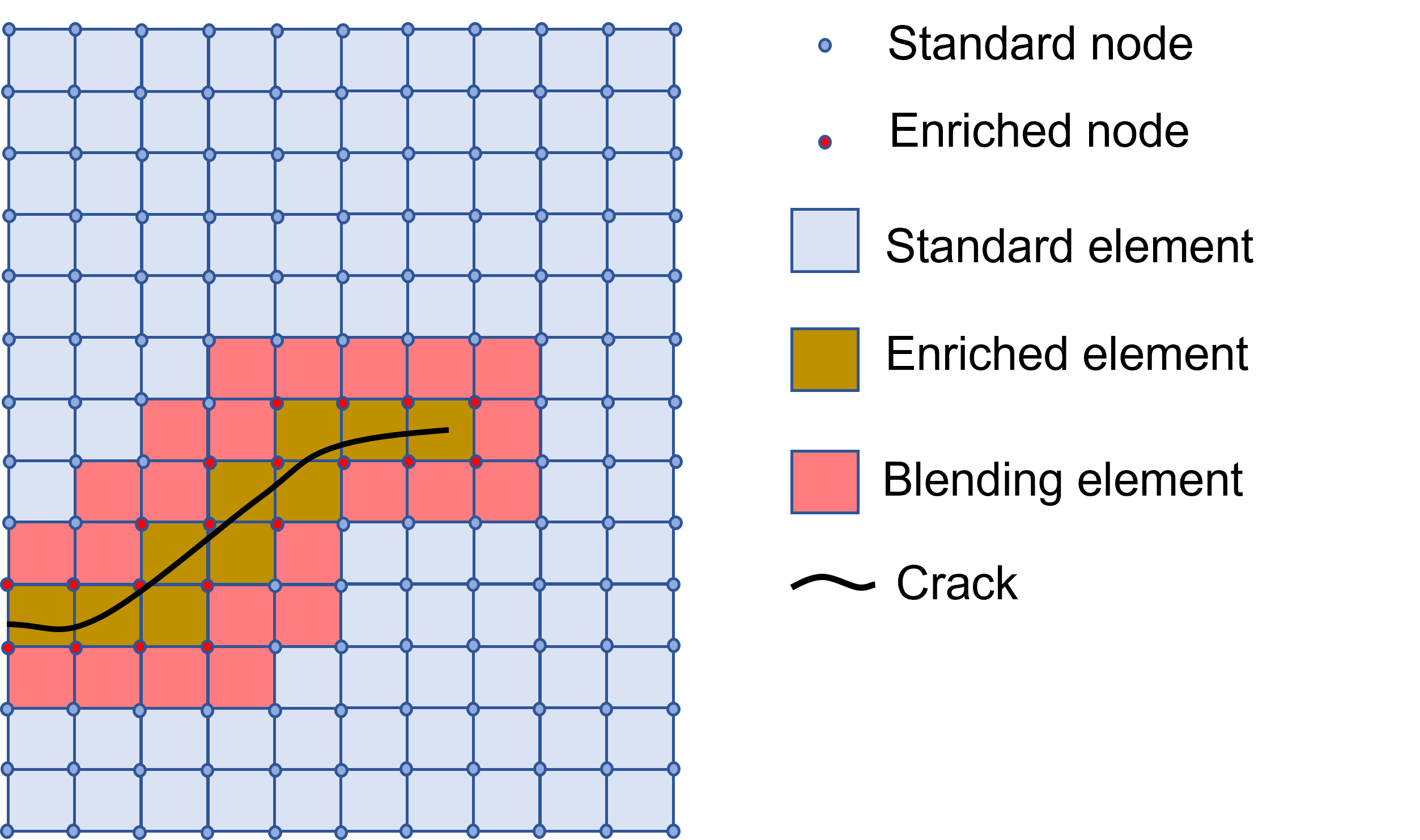}
    \caption{XFEM modeling of a cracked body. The Heaviside function enriches the elements across the crack line. Asymptotic functions enrich elements at the crack tips. Immediately neighboring elements of enriched elements, referred to as blending elements, are partially enriched.}
    \label{fig:xfemModel}
\end{figure}
 
 The extended finite element method (XFEM) offers a way of modeling discontinuities without remeshing. It was first introduced by \citet{belytschko_elastic_1999} and \citet{moes_finite_1999}, and later utilized in the commercial code of ABAQUS. Based on the concept of partition of unity proposed by \citet{melenk_partition_1996}, XFEM enriches the finite element formulation by the local enrichment function. The nodes of elements containing discontinuities are enhanced using the enrichment function with additional degrees of freedom (DOFs). Asymptotic functions extracted from analytical solutions in conjunction with the Heaviside jump function are utilized to model discontinuities across the crack line \cite{belytschko_arbitrary_2001}, as shown in \cref{fig:xfemModel}. All nodes in enriched elements, referred to as enriched nodes, are augmented with additional DOFs, while nodes in standard elements do not have such additional DOFs. The elements neighboring the enriched elements are referred to as blending elements, which are partially enriched since they have both enriched nodes and standard nodes.

Peridynamics (PD), introduced first by \citet{silling2000reformulation}, removes classical continuum theory's inadequacies in describing discontinuities. Instead of relying on derivatives of displacements, PD theory takes the force density in an integration form, making it suitable for analyzing structures containing discontinuities. However, the original version of the PD theory requires volume and surface correction to improve simulation accuracy. To remove these requirements, \citet{madenci_peridynamic_2019} developed a PD model based on least-squares minimization (PDLSM), which has been proved to recover the non-ordinary state based PD when the horizon is a sphere.  

Since PD is computationally expensive but can handle discontinuities, while FEM is well developed and computationally efficient, coupling PD with FEM is highly desirable for both methods.
\citet{sun_superposition-based_2019} devised a method of coupling based on a partial superposition of FEM and PD solutions for static and quasi-static problems.
\citet{shen_hybrid_2020} introduced truss elements to bridge finite element (FE) sub-regions and PD sub-regions to couple PD with FEM.
Using Lagrange multipliers, \citet{pagani_coupling_2020} developed a technique to couple 3-D peridynamics with 1-D high-order finite elements.  
Using the weighted residual method, the authors \cite{liu_revised_2021, liu_simulating_2021, liu_coupled_2021} proposed a straightforward framework to couple PD with FEM for 2-D and 3-D problems.
\citet{zaccariotto_coupling_2018, wang_yongwei_hybrid_2019} and \citet{tong_adaptive_2020} have developed various adaptive approaches which can transform FEM nodes into PD nodes for coupling PD with FEM.

In this paper, we extend the coupling framework proposed by the authors \cite{liu_simulating_2021, liu_coupled_2021} and present an adaptive coupling PDLSM-FEM for fracture analysis. The adaptive PDLSM-FEM is straightforward and does not require a transition zone to transform information between PD elements and finite elements. An adaptive algorithm is implemented with the dual goals of meshing regions with cracks using PD equations while minimizing the use of PD interactions to maximize computational efficiency. Within the presented framework, only elements intersecting with the crack and the neighboring elements are defined as PD elements whose stiffness matrices are derived based on the meshless method and PDLSM formulas, which is suitable for modeling the failure of structures. The remaining elements are conventional FEM elements. Numerical integration of interaction integral is proposed and implemented to evaluate the stress intensity factors (SIFs) for 2-D problems. The criterion of maximum circumferential tensile stress (MCTS) \cite{erdogan1963crack} is used for failure prediction. The crack growth direction $\theta_c$ is determined based on MCTS concerning the orientation, and crack starts to grow when the equivalent SIF $K_{eq}$ reaches the fracture toughness $K_{Ic}$. The PD elements are adaptively updated based on the current crack configuration and crack propagation. Although the criterion of MCTS used in this work can not handle crack branching and coalescence, it is better for quasi-static simulation of simple crack propagation in PD, compared to the bond stretch criterion \cite{silling_meshfree_2005} and bond energy criterion \cite{foster_energy_2011}, which may exist spurious bond breakage. Besides, the presented model does not require additional DOFs compared with the XFEM.

This paper is organized as follows. First, PDLSM theory is briefly reviewed in \cref{sec:pdlsm}, and the adaptive PDLSM-FEM is proposed and presented in \cref{sec:efempd}. After that, the 2-D simulation model of quasi-static crack growth is described in \cref{sec:scps}. Next, four examples are performed to demonstrate the proposed adaptive PDLSM-FEM in \cref{sec:NumR}. Finally, the
conclusions are drawn in \cref{sec:conc}.

\section{Peridynamics least-square minimization theory}\label{sec:pdlsm}

Since peridynamics was first introduced by \citet{silling2000reformulation}, various modifications and improvements to the theory have been proposed in the literature. The PDLSM introduced by \citet{madenci_peridynamic_2019} has some marked advantages, including allowing arbitrarily shaped interaction domains and eliminating surface correction. In this work, PDLSM is utilized to derive the element stiffness matrices for regions containing cracks. Detailed PDLSM can be found in \cite{madenci_peridynamic_2019}, but this section briefly reviews the formulation for completeness.

\begin{figure}[htbp]
    \centering
    \includegraphics[width=4in]{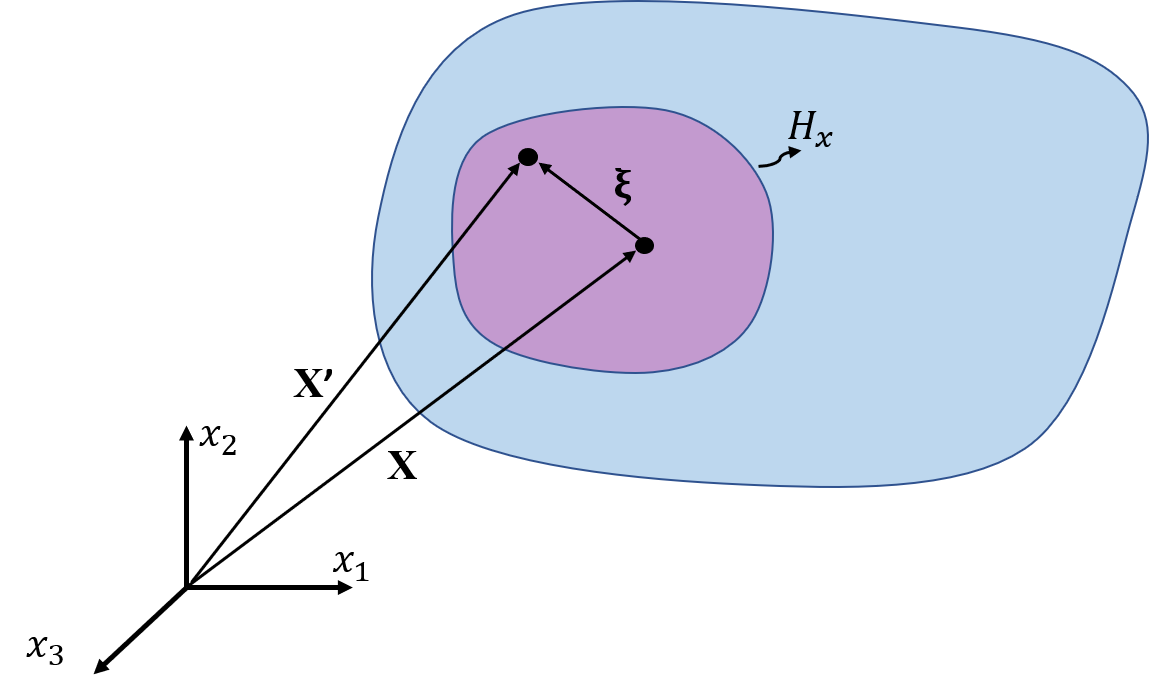}
    \caption{Field point $\textbf{x}$ interacts with its neighbors $\textbf{x}'$ within the interaction domain $H_x$.}
    \label{fig:PD_inter_doma}
\end{figure}

In PD,  point $\textbf{x}$ interacts with its neighbors $\textbf{x}'= \textbf{x}+\boldsymbol{\xi}$ within its interaction domain $H_x$, as illustrated in \cref{fig:PD_inter_doma}. Based on this nonlocal interaction concept of PD, the 2nd-order Taylor Series Expansion and the least-squares minimization \cite{madenci_peridynamic_2019} explicitly derived the nonlocal PDLSM differential operator as follows,
\begin{equation}\label{eq:diffOper}
    \left[ 
    \begin{matrix}
    \frac{\partial}{\partial x_1} & \frac{\partial}{\partial x_2} & \frac{\partial}{\partial x_3} & \frac{\partial^2}{\partial x_1^2} & \frac{\partial^2}{\partial x_2^2} & \frac{\partial^2}{\partial x_3^2} &\frac{\partial^2}{\partial x_1x_2} &\frac{\partial^2}{\partial x_2x_3} &\frac{\partial^2}{\partial x_3x_1}
    \end{matrix}
\right ]^T  f(\textbf{x}) =
\int_{H_x} \omega(\left|\boldsymbol{\xi}\right|) \left[ 
    \begin{matrix}
    \textbf{g} \\
    \textbf{d}
    \end{matrix}
\right]\left(f(\textbf{x}')-f(\textbf{x}) \right)dV'.
\end{equation}
The vectors $\textbf{g} = [ g_1 \quad g_2 \quad g_3 ]^T$ and  $\textbf{d} = [ d_1 \quad d_2 \quad d_3 \quad d_4 \quad d_5 \quad d_6 ]^T$ are defined as follows,
\begin{equation}\label{eq:gd}
    \left[
    \begin{array}{c}
         \textbf{g}  \\
         \textbf{d}
    \end{array}
    \right]
    =\textbf{A}^{-1} \hat{\boldsymbol{\xi}},
\end{equation}
where
\begin{equation}
    \hat{\boldsymbol{\xi}}=[\xi_1 \quad \xi_2 \quad \xi_3  \quad \xi_1^2 \quad \xi_2^2 \quad \xi_3^2 \quad \xi_1\xi_2 \quad \xi_2\xi_3 \quad \xi_1\xi_3]^T,
\end{equation}
\begin{equation}
    \textbf{A}=\left[
    \begin{matrix}
         \textbf{A}_{11}& \textbf{A}_{12} \\
         \textbf{A}_{21}& \textbf{A}_{22}
    \end{matrix}
    \right],
\end{equation}
\begin{equation}
    \textbf{A}_{11}=\int_{H_x} \omega(\left| \boldsymbol{\xi}\right|) \left[
    \begin{matrix}
        \xi_1^2     &  \xi_1\xi_2  & \xi_1\xi_3\\
         \xi_1\xi_2 &  \xi_2^2     & \xi_2\xi_3\\
         \xi_1\xi_3 &   \xi_2\xi_3 & \xi_3^2
    \end{matrix}
    \right] dV_{x'},
\end{equation}
\begin{equation}
    \textbf{A}_{12}=\int_{H_x} \omega(\left| \boldsymbol{\xi}\right|) \left[
    \begin{matrix}
       \frac{ \xi_1^3 }{2}    & \frac{\xi_1\xi_2^2}{2}   & \frac{\xi_1\xi_3^2}{2} &\xi_1^2\xi_2  &  \xi_1^2\xi_3 & \xi_1\xi_2\xi_3\\
        \frac{\xi_1^2\xi_2}{2}& \frac{ \xi_2^3 }{2}     & \frac{\xi_2\xi_3^2}{2} & \xi_1\xi_2^2 & \xi_1\xi_2\xi_3 & \xi_2^2\xi_3\\
         \frac{\xi_1^2\xi_3}{2} &  \frac{\xi_2^2\xi_3}{2}& \frac{ \xi_3^3 }{2} & \xi_1\xi_2\xi_3 & \xi_1\xi_3^2 & \xi_2\xi_3^2
    \end{matrix}
    \right] dV_{x'},
\end{equation}
\begin{equation}
    \textbf{A}_{21}=\int_{H_x} \omega(\left| \boldsymbol{\xi}\right|) \left[
    \begin{matrix}
        \xi_1^3 &  \xi_1^2\xi_2  & \xi_1^2\xi_3\\
         \xi_1\xi_2^2 &  \xi_2^3     & \xi_2^2\xi_3\\
         \xi_1\xi_3^2 &   \xi_2\xi_3^2 & \xi_3^3\\
         \xi_1^2\xi_2 & \xi_1\xi_2^2 & \xi_1\xi_2\xi_3\\
         \xi_1\xi_2\xi_3 & \xi_2^2\xi_3 & \xi_2\xi_3^2\\
         \xi_1^2\xi_3 & \xi_1\xi_2\xi_3 & \xi_1\xi_3^2
    \end{matrix}
    \right] dV_{x'},
\end{equation}
\begin{equation}
    \textbf{A}_{22}=\int_{H_x} \omega(\left| \boldsymbol{\xi}\right|) \left[
    \begin{matrix}
        \frac{\xi_1^4}{2}     &  \frac{\xi_1^2\xi_2^2}{2}  & \frac{\xi_1^2\xi_3^2}{2} & \xi_1^3\xi_2 & \xi_1^3\xi_3 & \xi_1^2\xi_2\xi_3 \\
         \frac{\xi_1^2\xi_2^2}{2} &  \frac{\xi_2^4}{2}    &\frac{\xi_2^2\xi_3^2}{2} & \xi_1\xi_2^3 & \xi_1\xi_2^2\xi_3 & \xi_2^3\xi_3\\
         \frac{\xi_1^2\xi_3^2}{2} &   \frac{\xi_2^2\xi_3^2}{2}&  \frac{\xi_3^4}{2} & \xi_1\xi_2\xi_3^2 &\xi_1\xi_3^3 & \xi_2 \xi_3^3\\
          \frac{\xi_1^3\xi_2}{2} & \frac{\xi_1\xi_2^3}{2} & \frac{\xi_1\xi_2\xi_3^2}{2} & \xi_1^2\xi_2^2 & \xi_1^2\xi_2\xi_3 & \xi_1\xi_2^2\xi_3\\
          \frac{\xi_1^2\xi_2\xi_3}{2} & \frac{\xi_2^3\xi_3}{2} & \frac{\xi_2\xi_3^3}{2} & \xi_1\xi_2^2\xi_3 & \xi_1\xi_2\xi_3^2 & \xi_2^2\xi_3^2\\
          \frac{\xi_1^3\xi_3}{2} &\frac{\xi_1\xi_2^2\xi_3}{2} &\frac{\xi_1\xi_3^3}{2}  &\xi_1^2\xi_2\xi_3 & \xi_1^2\xi_3^2 & \xi_1\xi_2\xi_3^2
    \end{matrix}
    \right] dV_{x'}.
\end{equation}

The nonlocal differential operator of \cref{eq:diffOper} was employed to derive the equation of motion, as pointed out in \cite{madenci_weak_2018}.
In PD theory, the equation of motion is defined as
\begin{equation}
    \rho \Ddot{\textbf{u}}=\textbf{L}^{pd}(\textbf{x},t)+\textbf{b}(\textbf{x},t),
\end{equation}
where $\textbf{b}$ represents the external body force, $\Ddot{\textbf{u}}$ is the acceleration, $\textbf{L}^{pd}$ represents the internal body force, and $\rho$ is the mass density. The internal body force $\textbf{L}^{pd}(\textbf{x},t)$ is calculated by the summation of bond force density between points $\textbf{x}$ and $\textbf{x}'$ over the domain $H_x$ of point $\textbf{x}$, as illustrated in \cref{fig:PD_inter_doma}. 
Based on the classical continuum mechanics, replacing $f(\textbf{x})$ in \cref{eq:diffOper} by displacement $u_i(\textbf{x})$, and introducing the bond status parameter $\mu$ gives the nonlocal displacement gradient,
\begin{equation}\label{eq:disp_grad}
    \nabla \textbf{u}^{pd}=\int_{H_x} \mu\omega(\left| \boldsymbol{\xi}\right|) \boldsymbol{\eta} \otimes \textbf{g} dV_{x'},
\end{equation}
and the nonlocal internal force vector,
\begin{equation}\label{eq:Lpd}
    \textbf{L}^{pd}=\nabla \cdot \boldsymbol{\sigma}^{pd}=\int_{H_x} \mu \omega(\left| \boldsymbol{\xi}\right|) \textbf{G} \boldsymbol{\eta} dV_{x'},
\end{equation}
where $\omega(\left| \boldsymbol{\xi}\right|)$ is the weight function, $\boldsymbol{\xi}=\textbf{x}'-\textbf{x}$ is the relative position,  $\boldsymbol{\eta}= \textbf{u}(\textbf{x}')-\textbf{u}(\textbf{x})$ is relative displacement, the bond status parameter $\mu$ is defined as
\begin{equation}
    \mu=
    \begin{cases}
         &1, \qquad \text{unbroken bond},\\
         &0, \qquad \text{broken bond},
    \end{cases}
\end{equation}
and the matrix $\textbf{G}$ is defined as below,
\begin{equation}\label{eq:G}
    \textbf{G}=\left[
    \begin{array}{ccc}
         (\lambda+\mathcal{G})d_1 +\mathcal{G} (d_1+d_2+d_3)& (\lambda+\mathcal{G})d_4 & (\lambda+\mathcal{G})d_6 \\
         (\lambda+\mathcal{G})d_4 & (\lambda+\mathcal{G})d_2 +\mathcal{G} (d_1+d_2+d_3) & (\lambda+\mathcal{G})d_5 \\
         (\lambda+\mathcal{G})d_6 & (\lambda+\mathcal{G})d_5 & (\lambda+\mathcal{G})d_3 +\mathcal{G} (d_1+d_2+d_3) 
    \end{array}
    \right],
\end{equation}

For 2-D problems, the vectors $\textbf{d}$ and $\textbf{g}$ are reduced to as $\textbf{d} = [ d_1 \quad d_2 \quad d_3]^T$ and $\textbf{g} = [ g_1 \quad g_2]^T$, respectively, and all the components with $\xi_3$ in vector $\hat{\boldsymbol{\xi}}$ and matrix $\textbf{A}$  will be removed. For plane stress cases, the matrix $\textbf{G}$ is reduced to as 
\begin{equation}
    \textbf{G}=\left[
    \begin{array}{ccc}
         \frac{E}{2(1-\nu)} d_1 +\mathcal{G} (d_1+d_2)& \frac{E}{2(1-\nu)}d_3  \\
         \frac{E}{2(1-\nu)}d_3 & \frac{E}{2(1-\nu)}d_2 +\mathcal{G} (d_1+d_2) 
    \end{array}
    \right],
\end{equation}
and for plane strain problems, 
\begin{equation}
    \textbf{G}=\left[
    \begin{array}{ccc}
         (\lambda+\mathcal{G})d_1 +\mathcal{G} (d_1+d_2)& (\lambda+\mathcal{G})d_3  \\
         (\lambda+\mathcal{G})d_3 & (\lambda+\mathcal{G})d_2 +\mathcal{G} (d_1+d_2) 
    \end{array}
    \right],
\end{equation}
here $\nu$ is the Poisson's ratio, $E$ is the Young's modulus, $\mathcal{G}$ is the shear modulus, and $\lambda$ is the Lame’s constant.

A critical, advantageous feature of PDLSM is that during the derivation of \cref{eq:disp_grad} and \cref{eq:Lpd}, the interaction domain $H_x$ is not required to be a sphere. As a result, the interaction domain for points near or on the boundary can be tailored to stay entirely inside the problem domain without causing any surface effect, and the volume correction, which is typically needed for PD models, is not required for PDLSM.  Another advantage is that PDLSM described the stress divergence $\nabla \cdot \boldsymbol{\sigma}$ in a nonlocal integral form (see \cref{eq:Lpd}) so that it does not require that points $\textbf{x}$ and $\textbf{x}'$ be part of their families as a pairwise bond.

\section{Adaptive PDLSM-FEM}\label{sec:efempd}

In this section, we propose the adaptive PDLSM-FEM in which PDLSM equations are used to model the regions containing cracks for discontinuity analysis.

\begin{figure}[htbp]
    \centering
    \includegraphics[width=4.5in]{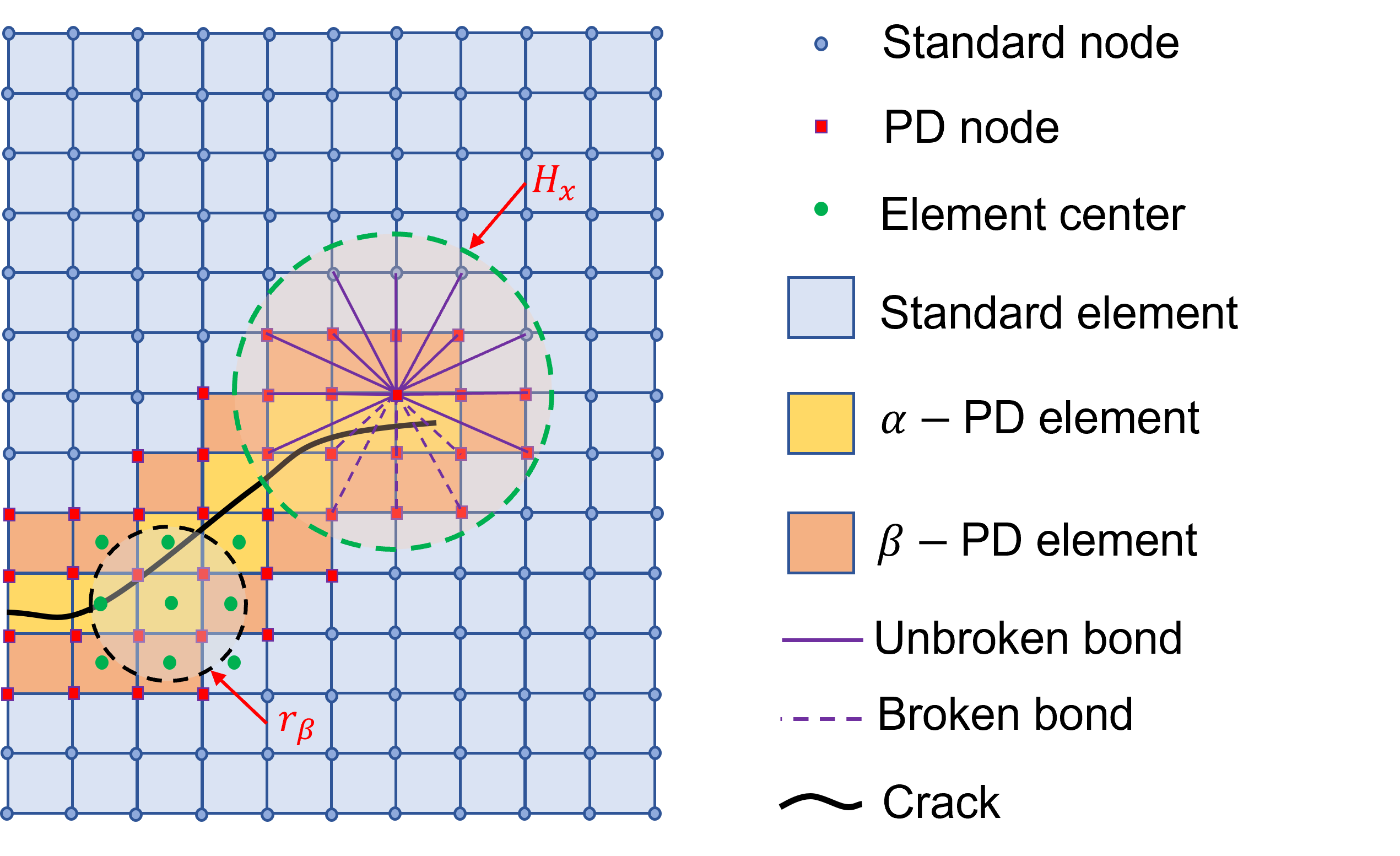}
    \caption{Adaptive PDLSM-FEM for modeling of a cracked body. The elements touched by the crack line are $\alpha$-PD elements. The $\beta$-PD elements' centers are located within circle $r_{\beta}$ centered at $\alpha$-PD elements' centers. The node of the PD element is a PD node that interacts with its neighbors within its interaction domain $H_x$.}
    \label{fig:efempd}
\end{figure}

In the framework of the adaptive PDLSM-FEM, the problem domain is meshed into elements, as presented in \cref{fig:efempd}.
 The crack line is explicitly described by line elements.
The elements touched by the crack line are referred to as $\alpha$-PD elements.
The $\beta$-PD elements are the elements whose centers are within the circles or spheres $r_{\beta}$ centered at $\alpha$-PD elements' centers.
We check the intersection between standard elements’ edges and crack line to determine the PD element; that is, if any edges of a standard element intersect with the crack line, it will become a $\alpha$-PD element, and its neighboring standard elements within the circle $r_\beta$ centered at its center will become $\beta$-PD elements.
The radius $r_{\beta}$ is $r_{\beta}=m_{\beta} \Delta_{min}$, where $\Delta_{min}$ is the minimum element size and $m_{\beta}$ is a constant.
$m_{\beta}<1$ means there are no $\beta$-PD elements. 

The nodes of PD elements are referred to as PD nodes governed by PDLSM equations. As presented in \cref{fig:efempd}, the PD node $\textbf{x}$ interacts with its family members, which are located within its interaction domain $H_x$. Whenever the interaction bond crosses the crack, it will be broken permanently. The remaining elements and nodes are referred to as standard elements and standard nodes, respectively. The element type depends on the proximity of the element to the current crack location and will be updated adaptively during crack propagation. 

The conventional FEM is employed within the standard elements, while within PD elements, the one-point quadrature for the PD equations is employed. To derive the resultant governing equations, the weighted residual method (WRM) is applied to the Neumann boundary conditions and the equations of motion over the problem body, as shown below,
\begin{equation}\label{eq:bal_WR}
    \int_{S_{tot}}\delta \textbf{u}^T(\boldsymbol{\sigma}\textbf{n}-\textbf{T})dS-\int_{V_{tot}} \delta\textbf{u}^T(\textbf{L}+\textbf{b}-\rho \Ddot{\textbf{u}})dV=0,
\end{equation}
where $\textbf{L}$ represents the internal body force, $\textbf{T}$ is the external traction, $S_{tot}$ and $V_{tot}$ represent the surface and volume of the whole body, respectively, $\textbf{n}$ represents the boundary's unit normal vector, and $\delta\textbf{u}$ represents the admissible displacement. The summation of elements' virtual inner work can express that of the domain as
\begin{equation}
     \int_{S_{tot}}\delta \textbf{u}^T\boldsymbol{\sigma}\textbf{n}dS-\int_{V_{tot}} \delta\textbf{u}^T \textbf{L}dV =\sum_{e=1}^{N} \delta U_e,
\end{equation}
where $N$ is the number of elements in the domain, and $\delta U_e$ is the virtual inner work of element and is defined as
\begin{equation}\label{eq:Ue}
    \delta U_e= \int_{S_e}\delta \textbf{u}^T\boldsymbol{\sigma}\textbf{n}dS-\int_{V_{e}} \delta\textbf{u}^T\textbf{L}dV,
\end{equation}
here $S_e$ is the element's surface and $V_e$ represents the element's volume.

\subsection{Stiffness matrix of standard element}\label{sec:Kestd}

The global element stiffness matrix $\textbf{K}_E^{std}$ of standard element based on conventional FEM formulas is presented in this section. As the conventional FEM is a well-established technology, only necessary formulas are presented in this section for completeness and without derivation for conciseness. 

Within each standard element, the internal force is defined in the differential form of
\begin{equation}
    \textbf{L}= \nabla \cdot \boldsymbol{\sigma}.
\end{equation}
Applying Gauss divergence theorem to \cref{eq:Ue} leads to the virtual inner work as
\begin{equation}\label{eq:dUstd}
     \delta U_e^{std}= \int_{S_e}\delta \textbf{u}^T\boldsymbol{\sigma}\textbf{n}dS-\int_{V_{e}} \delta\textbf{u}^T\textbf{L}dV=\int_{V_{e}}\nabla \delta \textbf{u}:\boldsymbol{\sigma} dV.
\end{equation}
Based on the conventional FEM formulas, the discrete form of the virtual work of standard elements is expressed as
\begin{equation}
    \delta U_e^{std}=\int_{V_{e}}\nabla \delta \textbf{u}:\boldsymbol{\sigma} dV=\delta \textbf{u}^T_e \textbf{K}^{std}_e \textbf{u}_e,
\end{equation}
here $\textbf{u}_e$ is the element's nodal displacement vector, $\textbf{K}^{std}_e$ is the local element stiffness matrix of the standard element and is defined as
\begin{equation}
    \textbf{K}^{std}_e=\int_{V_{e}} \textbf{B}^T \textbf{D} \textbf{B}dV,
\end{equation}
in which $\textbf{D}$ is elasticity matrix and is defined as follows. For 3-D problems, 
\begin{equation}
    \textbf{D}=\frac{E}{(1+\nu)(1-2\nu)}
    \left[
        \begin{matrix}
            (1-\nu) & \nu & \nu & 0 & 0 & 0\\
           \nu  & (1-\nu) & \nu & 0 & 0 & 0\\
            \nu & \nu &(1-\nu)  & 0 & 0 & 0\\
            0 & 0 & 0 & \frac{1-2\nu}{2} & 0 & 0\\
            0 & 0 & 0 & 0 & \frac{1-2\nu}{2} & 0\\
            0 & 0 & 0 & 0 & 0 & \frac{1-2\nu}{2}\\
        \end{matrix}
    \right],
\end{equation}
 For 2-D plane stress problem, 
\begin{equation}
    \textbf{D}=\frac{E}{1-\nu^2}
    \left[
        \begin{array}{ccc}
            1 & \nu & 0\\
           \nu  & 1 & 0\\
            0 & 0 &  \frac{1-\nu}{2}
        \end{array}
    \right],
\end{equation}
and for 2-D plane strain problem,
\begin{equation}
    \textbf{D}=\frac{E}{(1-2\nu)(1+\nu)}
    \left[
        \begin{array}{ccc}
            1-\nu & \nu & 0\\
           \nu  & 1-\nu & 0\\
            0 & 0 &  \frac{1-2\nu}{2}
        \end{array}
    \right].
\end{equation}
The matrix $\textbf{B}$ links the strains to the nodal displacements and is defined as
\begin{equation}
    \textbf{B}=\left[ \begin{matrix}
       \textbf{B}_1 & \textbf{B}_2 &\cdots &\textbf{B}_{N_e}
    \end{matrix}  \right]
\end{equation}
where $N_e$ is the number of element's nodes, and $\textbf{B}_i$ ($i=1,\ 2,\ \cdots,\  N_e$) is defined as follows. For 2-D problems,
\begin{equation}
    \textbf{B}_i=\left[ \begin{matrix}
       \frac{\partial \mathcal{N}_i}{\partial x_1} & 0 & \frac{\partial \mathcal{N}_i}{\partial x_2}\\
       0 & \frac{\partial \mathcal{N}_i}{\partial x_2} & \frac{\partial \mathcal{N}_i}{\partial x_1}
    \end{matrix}  \right]^T, 
\end{equation}
and for 3-D problems,
\begin{equation}
    \textbf{B}_i=\left[ \begin{matrix}
       \frac{\partial \mathcal{N}_i}{\partial x_1} & 0 & 0 & \frac{\partial \mathcal{N}_i}{\partial x_2} & 0 & \frac{\partial \mathcal{N}_i}{\partial x_3}\\
       0 & \frac{\partial \mathcal{N}_i}{\partial x_2} & 0 & \frac{\partial \mathcal{N}_i}{\partial x_1} &  \frac{\partial \mathcal{N}_i}{\partial x_3} & 0  \\
      0 & 0 & \frac{\partial \mathcal{N}_i}{\partial x_3} & 0 & \frac{\partial \mathcal{N}_i}{\partial x_2} & \frac{\partial \mathcal{N}_i}{\partial x_1}
    \end{matrix}  \right]^T,  
\end{equation}
in which $\mathcal{N}_i$ is the shape function. 

The element's nodal displacement vector $\textbf{u}_e$ can be mapped from the global nodal displacement vector $\textbf{u}_g$ by a mapping matrix $\textbf{M}_e$ as
\begin{equation}\label{eq:mapue}
    \textbf{u}_e=\textbf{M}_e \textbf{u}_g.
\end{equation}
 where the mapping matrix $\textbf{M}_e$ is not related to the coordinates and is only determined by the indices of element's DOFs in the global displacement vector as follows: if the $k$-th DOF of $\textbf{u}_e$ and the $j$-th DOF of $\textbf{u}_g$ are the same, $M_{e(kj)}=1$, otherwise $M_{e(kj)}=0$. Therefore, the virtual work of standard elements is transformed into
\begin{equation}
    \delta U_e^{std} =\delta \textbf{u}^T_g {\textbf{M}}^T_e\textbf{K}^{std}_e  {\textbf{M}}_e\textbf{u}_g =\delta \textbf{u}^T_g \textbf{K}^{std}_E\textbf{u}_g,
\end{equation}
from which the global element stiffness matrix $\textbf{K}^{std}_E$ follows
\begin{equation}
    \textbf{K}^{std}_E=\textbf{M}^T_e\textbf{K}^{std}_e  \textbf{M}_e.
\end{equation}

\subsection{Stiffness matrix of PD element}\label{sec:Keer}

In this section, the stiffness matrix of the PD element is formulated based on the PDLSM. The discrete form of the PDLSM equations is briefly presented first, and the derivation details can be found in \cite{liu_simulating_2021,liu_coupled_2021}. Next, the stiffness matrix of the PD element is derived.

\subsubsection{Discrete form of PD equations}
Consider a PD node $\textbf{x}_{(i)}$ has nodal displacement as
\begin{equation}
    \textbf{u}_{(i)}=
    \begin{cases}
        \begin{split}
            &\left[{u}_{1(i)} \quad {u}_{2(i)}\right]^T, \qquad \text{2-D},\\
            &\left[{u}_{1(i)} \quad {u}_{2(i)} \quad {u}_{3(i)}\right]^T, \qquad \text{3-D}.
        \end{split}
    \end{cases}
\end{equation}
The PD node $\textbf{x}_{(i)}$ interacts with its family members which are located within its interaction domain $H_{x_{(i)}}$. The displacement of the family members of the PD node $\textbf{x}_{(i)}$ is defined as
\begin{equation}
 \textbf{u}^{(i)}_f=\left[
 \begin{matrix}
           \textbf{u}^T_{(i)} & \textbf{u}^T_{(2)} & \textbf{u}^T_{(m)} & \cdots & \textbf{u}^T_{(N_{(i)})}
        \end{matrix} 
        \right]^T.
\end{equation}
here, $N_{(i)}$ represents the number of family members, and subscript $(m)$ represents the $m$th family member $\textbf{x}_{(m)}$ of node $\textbf{x}_{(i)}$. Note that $\textbf{x}_{(i)}$ is the first member. The vectors $\textbf{u}^{(i)}_f$ and $\textbf{u}_{(i)}$ can be mapped from the global nodal displacement vector $\textbf{u}_g$ through two mapping matrices $\textbf{M}_{(i)}$ and $\textbf{M}^{(i)}$ as follows,
\begin{equation}\label{eq:u_fami_toug2D}
    \textbf{u}^{(i)}_f=\textbf{M}^{(i)}\textbf{u}_g,
\end{equation}
\begin{equation}\label{eq:u_node_toug2D}
    \textbf{u}_{(i)}=\textbf{M}_{(i)}\textbf{u}_g,
\end{equation}
 in which the mapping matrices are obtained similar to the mapping matrix $\textbf{M}_e$ in \cref{eq:mapue}. From Eq. (\ref{eq:disp_grad}), we derive the strain vector of the PD node $\textbf{x}_{(i)}$ in a discrete form based on the simple one-point quadrature method as
\begin{equation}\label{eq:disc_strain2D}
\{ \varepsilon \}^{pd}_{(i)}= \textbf{C}^{(i)} \textbf{u}^{(i)}_f,
\end{equation}
where the matrix $\textbf{C}^{(i)}$ is defined as follows. For 2-D problems,
\begin{equation}\label{eq:C2D}
    \textbf{C}^{(i)}=\left[
       \begin{matrix}
         -\sum\limits_{m=2}^{N_{(i)}} C_{1(im)}& 0 &  C_{1(i2)}& 0 &\cdots & C_{1(iN_{(i))})}& 0 \\
         0 & -\sum\limits_{m=2}^{N_{(i)}}  C_{2(im)} & 0 &  C_{2(i2)} &\cdots & 0 & C_{1(iN_{(i))})}  \\
         -\sum\limits_{m=2}^{N_{(i)}} C_{2(im)}& -\sum\limits_{m=2}^{N_{(i)}} C_{1(im)} &C_{2(i2)} &C_{1(i2)} &\cdots &C_{2(iN_{(i))})} &C_{1(iN_{(i))})}
    \end{matrix} \right].
\end{equation}
and for 3-D problems,
\begin{equation}\label{eq:C3D}
    \textbf{C}^{(i)}=\left[
       \begin{matrix}
         -\sum\limits_{m=2}^{N_{(i)}} C_{1(im)}& 0 & 0 & C_{1(i2)}& 0 & 0 & \cdots  & C_{1(iN_{(i)})}& 0 & 0 \\
         0 & -\sum\limits_{m=2}^{N_{(i)}} C_{2(im)} & 0 & 0 & C_{2(i2)}  & 0  & \cdots  & 0 & C_{2(iN_{(i)})}  & 0\\
         0 & 0 & -\sum\limits_{m=2}^{N_{(i)}} C_{3(im)} & 0 & 0  & C_{3(i2)} & \cdots & 0 & 0  & C_{3(iN_{(i)})}\\
          -\sum\limits_{m=2}^{N_{(i)}} C_{2(im)}& -\sum\limits_{m=2}^{N_{(i)}} C_{1(im)} & 0 & C_{2(i2)}& C_{1(i2)} & 0 & \cdots & C_{2(iN_{(i)})}& C_{1(iN_{(i)})} & 0\\
           0 & -\sum\limits_{m=2}^{N_{(i)}} C_{3(im)} & -\sum\limits_{m=2}^{N_{(i)}} C_{2(im)} & 0 & C_{3(i2)}  & C_{2(i2)} & \cdots & 0 & C_{3(iN_{(i)})}  & C_{2(iN_{(i)})} \\
            -\sum\limits_{m=2}^{N_{(i)}} C_{3(im)} & 0 & -\sum\limits_{m=2}^{N_{(i)}} C_{1(im)} & C_{3(i2)} & 0  & C_{1(i2)} & \cdots  & C_{3(iN_{(i)})} & 0  & C_{1(iN_{(i)})}
    \end{matrix} \right],
\end{equation}
where $C_{k(im)} (k=1,2,3;\ m=2,3,\cdots, N_{(i)})$ are defined as
\begin{equation}\label{eq:ckm}
    C_{k(im)}=\mu_{(im)}\omega_{(im)} g_{k(im)}V_{(m)},
\end{equation}
where $g_{k(im)}$ is defined in \cref{eq:gd}, $V_{(m)}$ is the volume of the node $\textbf{x}_{(m)}$, $\mu_{(im)}$ and $\omega_{(im)}$ represents the bond status and the weight function, respectively.
Based on the constitutive law of classical continuum mechanics, the stress vector of the PD node $\textbf{x}_{(i)}$ in discrete form is
\begin{equation}\label{eq:disc_stress}
    \{\sigma\}^{pd}_{(i)}=\textbf{D}\textbf{C}^{(i)}\textbf{u}^{(i)}_f=\textbf{D}\textbf{C}^{(i)}\textbf{M}^{(i)}\textbf{u}_g,
\end{equation}

From \cref{eq:Lpd}, we derive the internal force vector of the enrich node $\textbf{x}_{(i)}$ in discrete form as follows:
\begin{equation}\label{eq:disc_Lpd}
    \textbf{L}^{pd}_{(i)}=\textbf{H}^{(i)} \textbf{u}^{(i)}_f = \textbf{H}^{(i)} \textbf{M}^{(i)}\textbf{u}_g,
\end{equation}
where $\textbf{H}$ is defined as
\begin{equation}\label{eq:H}
     \textbf{H}^{(i)}=\left[
    \begin{matrix}
    -\sum\limits_{m=2}^{N_{(i)}} \mu_{(im)}\omega_{(im)}\textbf{G}_{(im)}V_{(m)} &  \mu_{(i2)}\omega_{(i2)}\textbf{G}_{(i2)}V_{(2)} & \cdots & \mu_{(iN_{(i)})}\omega_{(iN_{(i)})}\textbf{G}_{(iN_{(i)})}V_{(N_{(i)})}
    \end{matrix}
    \right].
\end{equation}

\subsubsection{Stiffness matrix of PD element}

Following \cref{eq:Ue}, the virtual inner work of PD elements is 
\begin{equation}\label{eq:dUer}
     \delta U_e^{pd}= \int_{S_e}\delta \textbf{u}^T\boldsymbol{\sigma}\textbf{n}dS-\int_{V_{e}} \delta\textbf{u}^T\textbf{L}dV.
\end{equation}

Within the PD element, we do not employ the Gauss divergence theorem and Gauss integration as we do with the standard element, since discontinuities may exist within the PD element for which the divergence theorem and Gauss integration are invalid. Instead, we apply the one-point quadrature method, which was employed in many published PD models, and the nonlocal stress (\cref{eq:disc_stress}) and nonlocal internal force (\cref{eq:disc_Lpd}) are used.

Based on \cref{eq:disc_Lpd} and the one-point quadrature method, the virtual work of the inner body force of PD element is 
\begin{equation}
    \begin{split}
        \delta U_{eB}^{pd}&=-\int_{V_{e}} \delta\textbf{u}^T\textbf{L}dV\\
        &=-\sum_{i=1}^{N_e} \frac{V_e}{N_e} \delta \textbf{u}_{(i)}^T \textbf{L}^{pd}_{(i)} \\
        &=-\sum_{i=1}^{N_e} \frac{V_e}{N_e} \delta \textbf{u}_g^T \textbf{M}_{(i)}^T \textbf{H}^{(i)} \textbf{M}^{(i)}\textbf{u}_g\\
        &=\delta \textbf{u}_g^T \left(-\sum_{i=1}^{N_e} \frac{V_e}{N_e} \textbf{M}_{(i)}^T \textbf{H}^{(i)} \textbf{M}^{(i)}\right) \textbf{u}_g\\
        &=\delta \textbf{u}_g^T \textbf{K}^{pd}_{EB} \textbf{u}_g,
    \end{split}
\end{equation}
in which $V_e$ represents the element volume, $N_e$ represents the number of nodes of element, $\textbf{K}^{pd}_{EB}$ is the stiffness matrix from the interior body force contribution defined as
\begin{equation}
    \textbf{K}^{pd}_{EB}=-\sum_{i=1}^{N_e} \frac{V_e}{N_e} \textbf{M}_{(i)}^T \textbf{H}^{(i)} \textbf{M}^{(i)}.
\end{equation}
The interior surface traction of the PD element can be expressed as
\begin{equation}\label{eq:tl_mat}
    \boldsymbol{\sigma}\textbf{n}=\textbf{N}\{ \sigma\}.
\end{equation}
For 2-D problems, the normal matrix $\textbf{N}$ is defined as
\begin{equation}\label{eq:N2d}
    \textbf{N}=\left[
    \begin{matrix}
    n_{1} & 0 & n_{2}\\
    0 & n_{2} & n_{1}
    \end{matrix}
    \right].
\end{equation}
 and for 3-D problems,
 \begin{equation}\label{eq:N3d}
    {\textbf{N}}=\left[
    \begin{matrix}
   {n}_{1} & 0 & 0 &  {n}_{2} & 0 &  {n}_{3}\\
    0 &  {n}_{2} & 0 &  {n}_{1} &  {n}_{3} & 0\\
    0 & 0 &  {n}_{3} & 0 &  {n}_{2} &  {n}_{1}
    \end{matrix}
    \right].
\end{equation}
The virtual work of the internal surface traction of the PD element is
\begin{equation}
    \begin{split}
        \delta U_{eS}^{pd} &= \int_{S_e}\delta \textbf{u}^T\boldsymbol{\sigma}\textbf{n}dS\\
        &=\sum^{N_s}_{s}\sum_i^{N_{sn}}\frac{A_s}{N_{sn}} \delta \textbf{u}^T_{(i)} \textbf{N}_s \{ \sigma\}^{pd}_{(i)} \\
        &=\sum^{N_s}_{s}\sum_i^{N_{sn}}\frac{A_s}{N_{sn}} \delta \textbf{u}^T_g \textbf{M}^T_{(i)} \textbf{N}_s \textbf{D}\textbf{C}^{(i)}\textbf{M}^{(i)}\textbf{u}_g\\
        &=\delta \textbf{u}^T_g \textbf{K}^{pd}_{ES} \textbf{u}_g
    \end{split}
\end{equation}
in which $N_s$ is the number of the surface of the element, $N_{sn}$ is the number of the nodes on the surface $s$, $A_s$ is the area of the surface $s$, $\textbf{N}_s$ is the normal matrix of the surface $s$ as defined in \cref{eq:N2d} or \cref{eq:N3d}, and $\textbf{K}^{pd}_{ES}$ is the stiffness matrix from surface traction contribution and is defined as
\begin{equation}
    \textbf{K}^{pd}_{ES}=\sum^{N_s}_{s}\sum_i^{N_{sn}}\frac{A_s}{N_{sn}}  \textbf{M}^T_{(i)} \textbf{N}_s \textbf{D}\textbf{C}^{(i)}\textbf{M}^{(i)}.
\end{equation}
Thus, the total virtual internal work within the PD element is
\begin{equation}
    \delta U_e^{pd}=\delta U_{eB}^{pd}+\delta U_{eS}^{pd}=\delta \textbf{u}^T_g \left( \textbf{K}^{pd}_{EB}+\textbf{K}^{pd}_{ES}\right) \textbf{u}_g =\delta \textbf{u}^T_g \textbf{K}^{pd}_E \textbf{u}_g,
\end{equation}
in which $\textbf{K}^{pd}_E$ is the global element stiffness of the PD elements defined as
\begin{equation}
    \textbf{K}^{pd}_E=\textbf{K}^{pd}_{EB}+\textbf{K}^{pd}_{ES}.
\end{equation}
It is worthy to note that if an element surface is shared by two PD elements, its surface traction has no contribution to the total virtual internal work, since the normal vector is $\textbf{n}$ on one side, and is $(-\textbf{n})$ on the other side. Therefore, the virtual work of both will cancel each other.

\subsection{Governing equations}

As presented in \cref{sec:Kestd,sec:Keer}, the total virtual internal  work is 
\begin{equation}\label{eq:Uinner}
\begin{split}
      \delta U &=\sum_e^{N^{std}} \delta U_e^{std} + \sum_e^{N^{pd}} \delta U_e^{pd}\\
      &= \delta \textbf{u}^T_g \left(\sum_e^{N^{std}} \textbf{K}^{std}_E \right) \textbf{u}_g + \delta \textbf{u}^T_g \left(\sum_e^{N^{pd}} \textbf{K}^{pd}_E \right) \textbf{u}_g\\
      &=\delta \textbf{u}^T_g \textbf{K} \textbf{u}_g
\end{split}
\end{equation}
in which $N^{std}$ is the total number of standard element, $N^{pd}$ is the total number of PD elements, and $\textbf{K}$ is the global stiffness matrix defined as
\begin{equation}
    \textbf{K}=\sum_e^{N^{std}} \textbf{K}^{std}_E + \sum_e^{N^{pd}} \textbf{K}^{pd}_E.
\end{equation}
Furthermore, the virtual work from the inertial term can be transformed as 
\begin{equation}\label{eq:Uitot}
    \delta U_I=\int_{V_{tot}} \delta\textbf{u}^T\rho \Ddot{\textbf{u}}dV=\delta \textbf{u}_g^T\textbf{M}\Ddot{\textbf{u}}_g,
\end{equation}
and the virtual work from the external forces $\textbf{b}$ and $\textbf{T}$ can be transformed as
\begin{equation}\label{eq:UF}
    \delta W=\int_{V_{tot}} \delta\textbf{u}^T\textbf{b}dV + \int_{S_{tot}} \delta \textbf{u}^T\textbf{T}dS=\delta \textbf{u}_g^T\textbf{F},
\end{equation}
Substituting \cref{eq:Uinner,eq:Uitot,eq:UF} into \cref{eq:bal_WR} leads to the governing equations of the adaptive PDLSM-FEM model as follows,
\begin{equation}\label{eq:gove_equa}
    \textbf{M}\Ddot{\textbf{u}}_g + \textbf{K}\textbf{u}_g=\textbf{F},
\end{equation}
For quasi-static and static problems, Eq. (\ref{eq:gove_equa}) is reduced to as 
\begin{equation}\label{eq:kuf}
    \textbf{K}\textbf{u}_g=\textbf{F}.
\end{equation}
In this study, the stiffness matrix \textbf{K} is assembled and stored in a compressed sparse row format with only storing the non-zero coefficients and their row and column information. The assembling algorithm can be found in an earlier publication by the authors \citep{liu_coupled_2021}.

\section{simulating 2-D quasi-static crack growth}\label{sec:scps}
This section presents the quasi-static crack growth simulation model for 2-D problems based on the adaptive PDLSM-FEM and the linear elastic fracture mechanics (LEFM). The crack propagation criterion is based on the MCTS evaluated from SIFs within the framework of LEFM, where the SIFs are calculated using the $I$-integral \citep{yau_mixed-mode_1980}, which is derived from the $J$-integral. The evaluation of SIFs using the $I$-integral is first discussed in \cref{sec:SIFs}. After that, the numerical integration of $I$-integral is described in \cref{sec:NumJ}. Next, the failure criterion is described in \cref{sec:spc}. Finally, the step-by-step  simulation process is presented in \cref{sec:simPro}.

\subsection{Evalution of SIFs}\label{sec:SIFs}

 In this work, the SIFs are evaluated based on the $I$-integral which is an extension of $J$-integral. The concept, first proposed by \citet{yau_mixed-mode_1980}, is reviewed briefly. 

Consider two states of a body with crack, that are, the current state (1) denoted by $\left(u_i^{(1)}, \varepsilon_{ij}^{(1)},\sigma_{ij}^{(1)} \right)$ and an auxiliary state (2) called by $\left(u_i^{(2)}, \varepsilon_{ij}^{(2)},\sigma_{ij}^{(2)} \right)$ , the $I$-integral is defined as
\begin{equation}\label{eq:ii}
    I^{(1,2)}=\int_{\Gamma} \left(w^{(1,2)}\mathcal{D}_{1j}-\sigma_{ij}^{(1)}\frac{\partial u_i^{(2)}}{\partial x_1} -\sigma_{ij}^{(2)}\frac{\partial u_i^{(1)}}{\partial x_1}\right)n_j d \Gamma,
\end{equation}
where $\mathcal{D}$ represents the Kronecker delta, $n_j$ represents the normal vector of the integration contour $\Gamma$, ($x_1,\ x_2$) are local coordinates centered at the crack tip, and $w^{(1,2)}$ is defined as below,
\begin{equation}
     w^{(1,2)}=\sigma_{ij}^{(1)}\varepsilon_{ij}^{(2)}=\sigma_{ij}^{(2)}\varepsilon_{ij}^{(1)}.
\end{equation}
Choosing state (2) as the asymptotic fields with $K_I^{(2)}=1$ and $K_{II}^{(2)}=0$ leads to the SIF $K_I^{(1)}$ for the current state (1),
\begin{equation}\label{eq:K1}
    K_I^{(1)}=\frac{E^*}{2} I^{(1,2)}.
\end{equation}
Similarly, letting the state (2) be the asymptotic fields with $K_I^{(2)}=0$ and $K_{II}^{(2)}=1$ leads to the SIF $K_{II}^{(1)}$ for the current state (1) as
\begin{equation}\label{eq:K2}
    K_{II}^{(1)}=\frac{E^*}{2} I^{(1,2)}.
\end{equation}
where $E^*$ is defined as
\begin{equation}
   E^*=
\begin{cases}
E,\quad \mathrm{plane\ stress}, \\
\frac{E}{1-\nu^2},\quad \mathrm{plane\ strain}.
\end{cases}
\end{equation}

\subsection{Numerical Integration of $I$-Integral}\label{sec:NumJ}

\begin{figure}[htbp]
    \centering
    \includegraphics[width=4in]{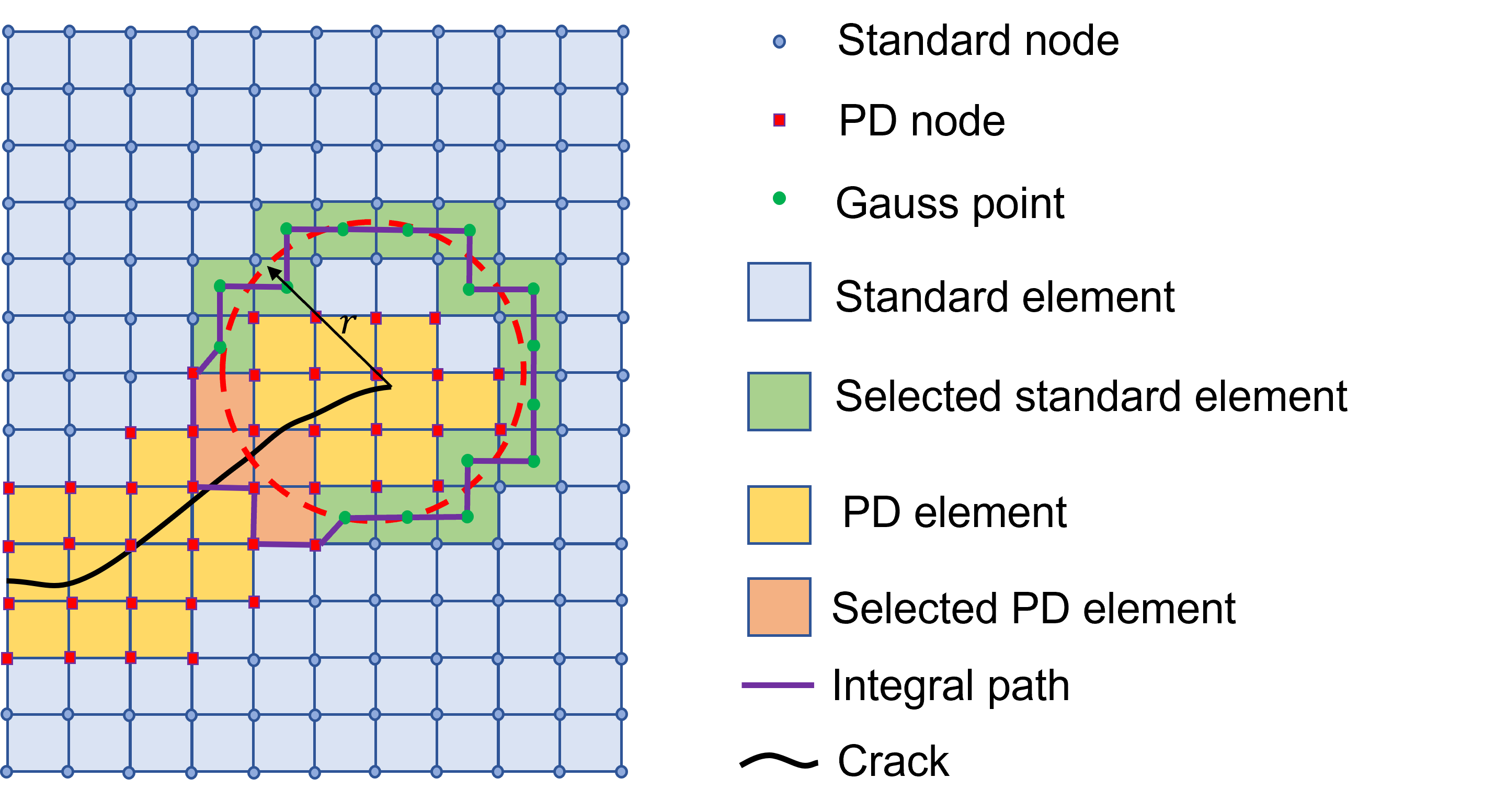}
    \caption{The contour of the $I$-integral. Elements intersected by the base circle (red dot line) are selected elements. Connecting consecutively all PD nodes of selected PD elements (light orange elements) that are outside of the base circle and the Gauss points of selected standard elements (light green element), constructs the integral contour (solid purple lines).}
    \label{fig:jPath}
\end{figure}

Although it is a common practice to evaluate the $J$- or $I$-integral by transforming the contour integral to an equivalent area integral and then using the Gauss integration over the elements of the equivalent area \citep{li_comparison_1985}, the approach may be invalid for elements cut through with cracks, as Gauss integration requires the integrand to be continuous. To evaluate SIFs, we develop in this section a new numerical method of $I$-integral over a contour of PD nodes and Gauss points of the standard element, based on the rationale that stresses and strains at Gauss points are more accurate than other locations for standard elements. Stresses and strains of PD nodes can be directly computed based on \cref{eq:disc_strain2D,eq:disc_stress}.

To determine the contour, we define a base circle with radius $r=m_r\Delta_{min}$ centered at the crack tip, as shown in \cref{fig:jPath}, in which $m_r$ is a constant and $\Delta_{min}$ is the minimum element size. Next, standard and PD elements intersected by the circle are designated as selected elements, as illustrated in \cref{fig:jPath}. The integral contour is constructed by connecting consecutively all PD nodes of selected PD elements (light orange elements) outside the base circle and the Gauss points of selected standard elements (light green elements), as presented by the solid purple lines in \cref{fig:jPath}. With the contour formed by the selected elements, we calculate the interaction integral $I^{(1,2)}$ as
\begin{equation}\label{eq:I12}
    I^{(1,2)}=\sum_{p=1}^{N_{p}}\left(F_{(p)}+F_{(p+1)} \right)\frac{l_p}{2},
\end{equation}
in which $l_p$ is the length between points $(p)$ and $(p+1)$, and $N_{p}$ represents the number of points on the integral contour. Point $(N_{p}+1)$ is identical to point $(1)$, ensuring the closed integral contour. $F_{(p)}$ denotes the integrand at point $(p)$ as:
\begin{equation}\label{eq:Fp}
    F_{(p)}=\left(w^{(1,2)}\mathcal{D}_{1j}-\sigma_{ij}^{(1)}\frac{\partial u_i^{(2)}}{\partial x_1} -\sigma_{ij}^{(2)}\frac{\partial u_i^{(1)}}{\partial x_1}\right)n_j.
\end{equation}

To evaluate the SIFs by \cref{eq:K1,eq:K2}, the displacement, strain, and stress fields $u_i^{(1)}$, $\varepsilon_{ij}^{(1)}$, and $\sigma_{ij}^{(1)}$ of state (1) are calculated from the adaptive PDLSM-FEM model, and $u_i^{(2)}$, $\varepsilon_{ij}^{(2)}$, and $\sigma_{ij}^{(2)}$ of the state (2) are calculated from LEFM formulas with $K_I^{(2)}=1$, $K_{II}^{(2)}=0$ or $K_I^{(2)}=0$, $K_{II}^{(2)}=1$. Then, the $I$-integral $I^{(1,2)}$ is calculated by \cref{eq:I12}. 

\begin{figure}[htbp]
    \centering
    \includegraphics[width=4in]{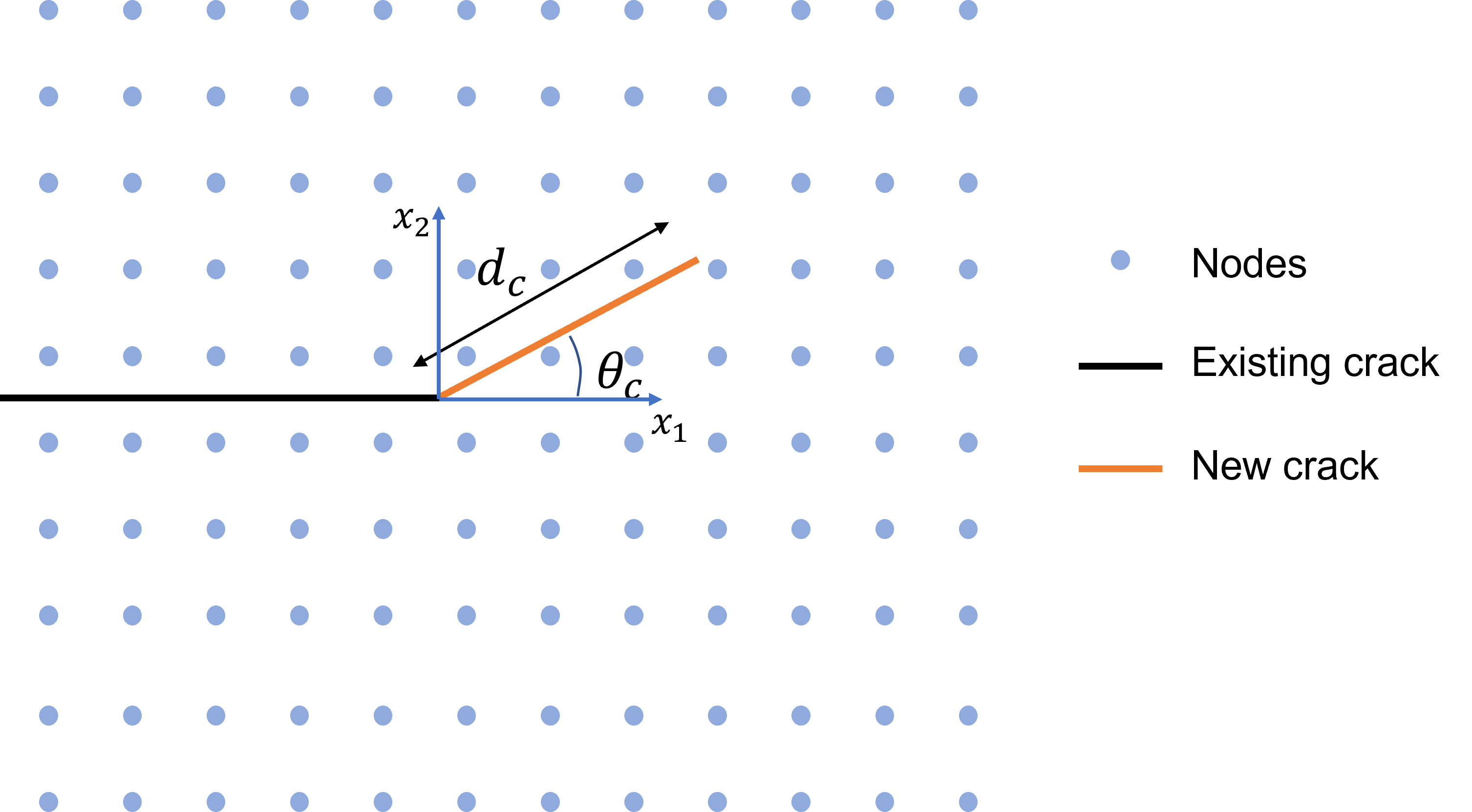}
    \caption{Schematic of crack growth. $\theta_c$ is the crack growth direction and $d_c$ is the crack growth amount at each step.}
    \label{fig:cracGrow}
\end{figure}

\subsection{Crack propagation criterion}\label{sec:spc}

To simulate quasi-static crack growth, three conditions must be determined: (1) the onset of growth, (2) the direction of growth, and (3) the amount of growth. 
\begin{figure}[htbp]
    \centering
    \includegraphics[width=4.8in]{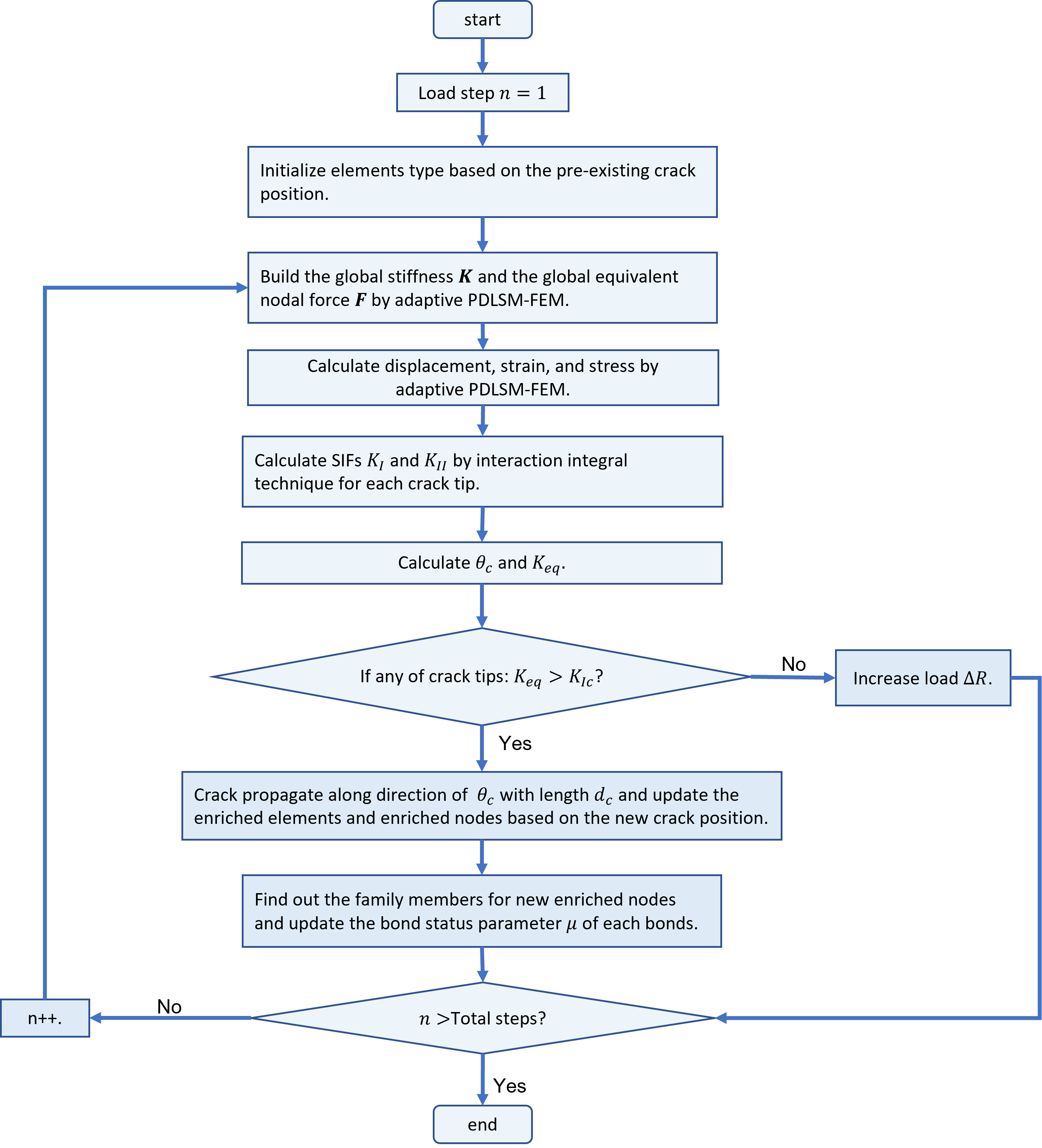}
   \caption{Flowchart of crack propagation simulation.}
    \label{fig:ProccracGrow}
\end{figure}

Various criteria have been developed for mixed-mode loading to determine the first two conditions. In the current work, the MCTS theory \citep{erdogan1963crack} is used to determine the onset and the direction of crack growth. For LEFM, the singular asymptotic stresses at the crack tip from LEFM solutions are used in the MCTS criterion to find out the crack growth direction $\theta_c$ (\cref{fig:cracGrow}), which is defined as
\begin{equation}\label{eq:thetaC}
\theta_c=
\begin{cases}
2\tan^{-1}{\left(\frac{K_I}{4K_{II}}-\frac{1}{4}\sqrt{\left(\frac{K_I}{K_{II}}\right)^2+8}\right)}\quad \mathrm{for} K_{II}>0, \\
2\tan^{-1}{\left(\frac{K_I}{4K_{II}}+\frac{1}{4}\sqrt{\left(\frac{K_I}{K_{II}}\right)^2+8}\right)}\quad \mathrm{for} K_{II}< 0,
\end{cases}
\end{equation}
For onset of crack growth, the MCTS must reach a critical value, which results in the equivalent SIF as
\begin{equation}\label{eq:keq}
    K_{eq}=K_I\cos^3{\frac{\theta_c}{2}}-\frac{3}{2}K_{II}\cos{\frac{\theta_c}{2}}\sin{\theta_c}.
\end{equation}
with $\theta_c$ determined by \cref{eq:thetaC}. If $K_{eq}>K_{Ic}$, the crack will grow in the direction of $\theta_c$ with an amount of $d_c$, as shown in \cref{fig:cracGrow}. In this work, $d_c$ is set to be $\alpha\Delta_{min}$, where $\Delta_{min}$ is the minimum element size, and $\alpha \in[0.5,1.5]$.

\subsection{Simulation process and flowchart}\label{sec:simPro} 

This section presents the quasi-static crack growth simulation step-by-step to obtain an overall sense of the adaptive PDLSM-FEM model. Fig. \ref{fig:ProccracGrow} shows the flowchart of the simulation process. 
\begin{itemize}
    \item Step 1. Initialize all elements as standard elements first. The elements intersecting with the pre-existing crack and their neighboring elements are transformed into PD elements. A small load is then applied to the material domain, and the nodal force vector $\textbf{F}$ and the stiffness matrix $\textbf{K}$ are assembled as described in \cref{sec:efempd}. 
    \item Step 2. The displacement field is solved by $\textbf{U}=\textbf{K}^{-1}\textbf{F}$, from which the stress and strain fields are then evaluated. The stresses and strains of the PD nodes are evaluated by PD \cref{eq:disc_strain2D,eq:disc_stress}, while strains and stresses in the standard elements are obtained based on conventional FEM. 
    \item Step 3. With the fields of displacements, strains, and stresses, the SIFs, $K_I$ and $K_{II}$ are computed from \cref{eq:K1,eq:K2,eq:I12,eq:Fp} for the current crack configuration, and then the $\theta_c$ and $K_{eq}$ are evaluated by \cref{eq:thetaC} and \cref{eq:keq}, respectively. 
    \item Step 4. If $K_{eq}<K_{Ic}$, the crack will propagate along with the direction $\theta_c$ by amount $d_c$. The PD elements are then updated based on the new crack location. All bonds across the new crack will be treated as broken. The global stiffness K is then updated, which becomes more compliant. 
    If $K_{eq}>K_{Ic}$, a small load will be added by
    $\Delta R= \text{max} \left\{\left(\frac{K_{Ic}}{K_{eq}}-1\right)R_n,\ \Delta R_{max}\right\}$, where $R_n$  is the current loading value, and $\Delta R_{max}$ is the pre-defined maximum loading increment.
    \item Step 5. The simulation goes to the next step. 
\end{itemize}

\section{Numerical results}\label{sec:NumR}

To validate the capability of the adaptive PDLSM-FEM model presented in this work, simulations of three 2-D plane stress problems and one 3-D problem are performed and presented in this section. In the first problem, an infinite plate with an inclined stationary crack is subjected to remote uniform traction leading to mixed-mode loading. In the second problem, a diagonal plate with an initial inclined crack is loaded under displacement control, and the crack growth is modeled using a quasi-static simulation approach. In the third problem, a compact tension test specimen with an initial crack in the plane of symmetry is loaded under displacement control. The crack growth is modeled again using the quasi-static approach. In the fourth problem, a 3-D block with a stationary crack is loaded under displacement control. The first, second, and fourth problems have been simulated and reported in an earlier publication using a general method of coupled PDLSM with FEM \citep{liu_simulating_2021,liu_coupled_2021}, and they are studied here again using the new adaptive PDLSM-FEM model. The simulations are performed using the PDLSM-FEM solver \cite{Liu2021_solver}.

 In this work, the weight function is specified in the Gauss distribution form as
\begin{equation}
    \omega(\left|\boldsymbol{\xi} \right|)=e^{-(\frac{\left|\boldsymbol{\xi} \right|}{c \delta_{(i)}})^2},
\end{equation}
here, $\delta_{(i)}$ is the horizon size of node $\textbf{x}_{(i)}$, $c$ is a constant and is specified as $c=1/3$ in this work based on experience.

\begin{figure}[htbp]
    \centering
         \subfloat[]{\label{fig:infiniPlatGemo}
    \includegraphics[width=2.in]{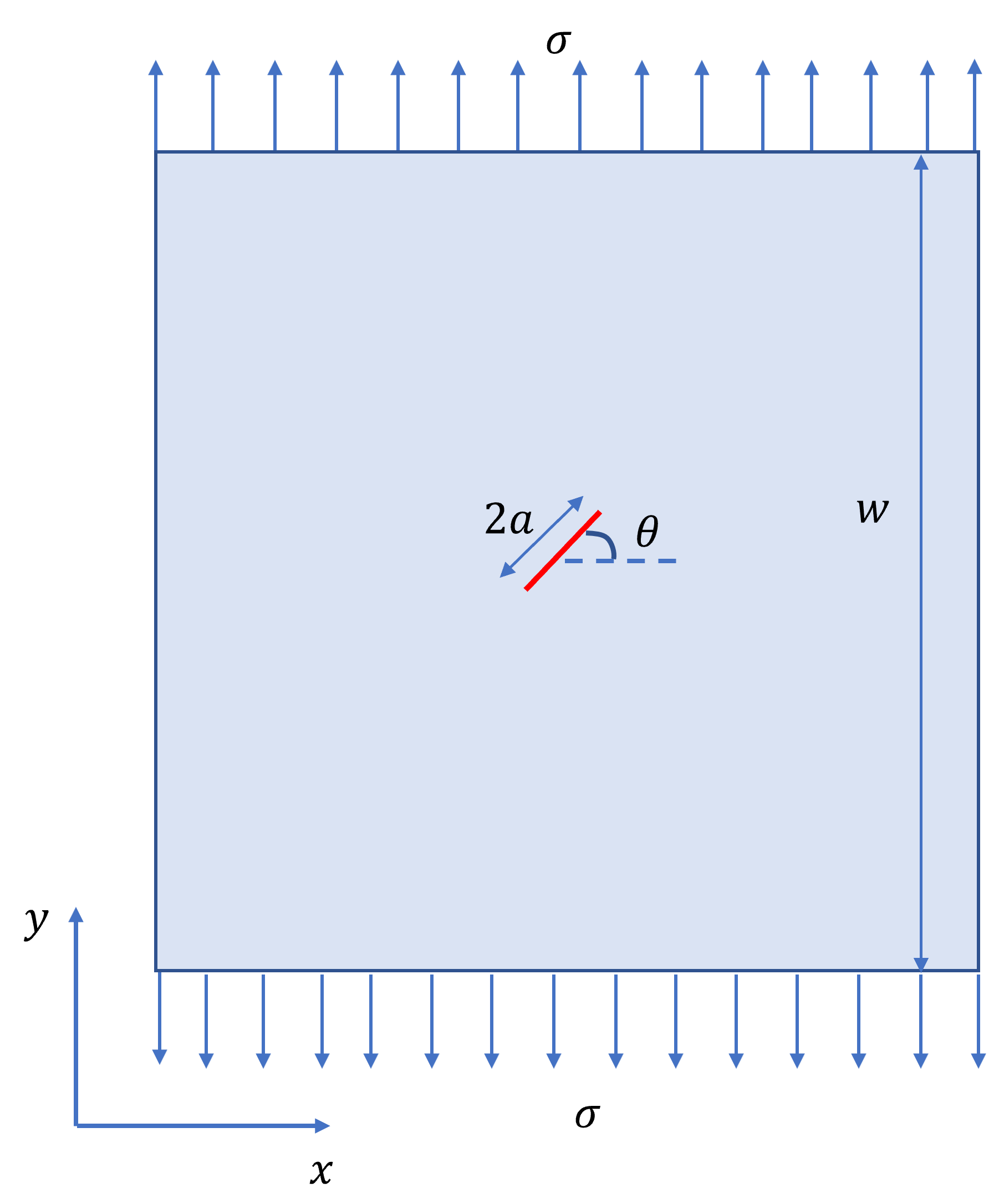}}\quad \subfloat[]{\label{fig:infiniPlatMSH}
    \includegraphics[width=2.4in]{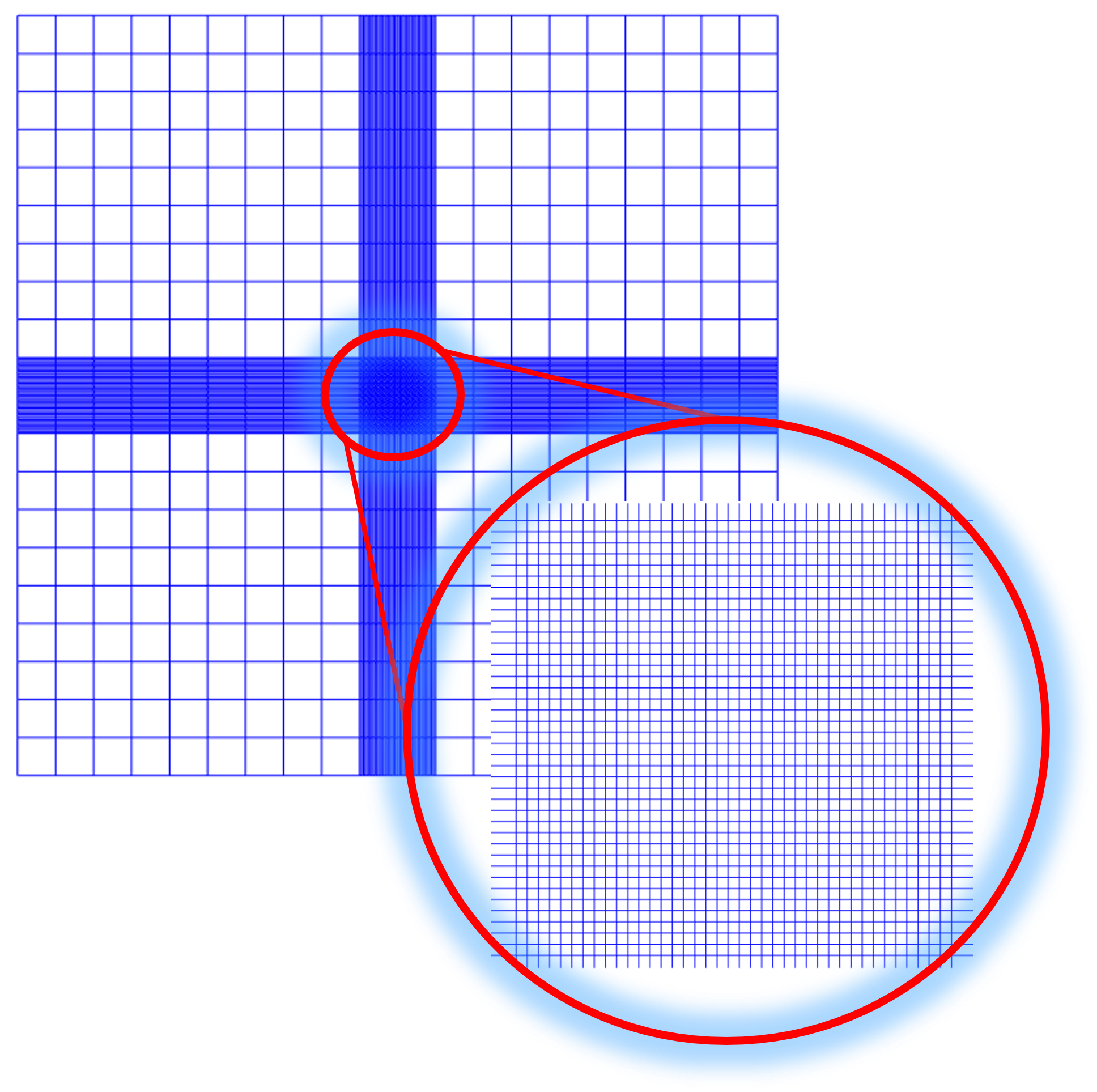}}
    \caption{A 2-D infinite plate under mixed-mode loading: (a) Geometry and boundary conditions, (b) Mesh.}
    \label{fig:infinPlat}
\end{figure}

\begin{figure}[htbp]
    \centering
    \includegraphics[width=3in]{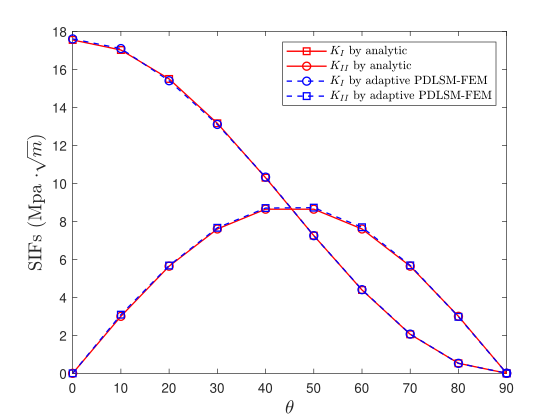}
    \caption{Numerical and analytical evaluation of SIFs at various crack inclined angles $\theta$.}
    \label{fig:SIFs}
\end{figure}
\begin{figure}[htbp]
    \centering
         \subfloat[]{\label{fig:infinitPlaSx}
    \includegraphics[width=1.5in]{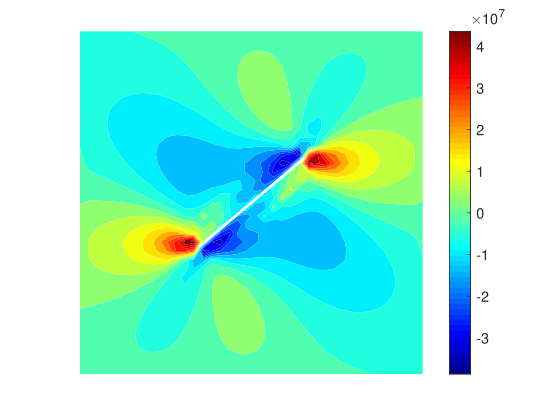}}\  \subfloat[]{\label{fig:infinitPlaSy}
    \includegraphics[width=1.5in]{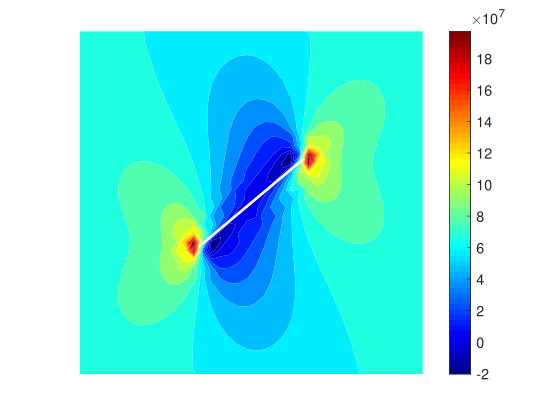}}\ 
    \subfloat[]{\label{fig:infinitPlaSxy}
    \includegraphics[width=1.5in]{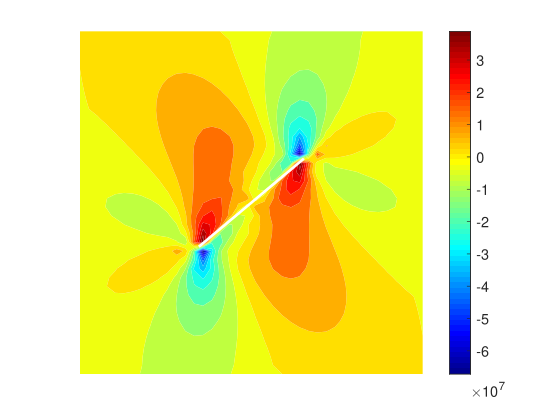}}
    \caption{Stress distributions around the crack of the 2-D infinite plate with crack angle $\theta=40^{\circ}$: (a) $\sigma_{xx}$, (b) $\sigma_{yy}$, (c) $\sigma_{xy}$.}
    \label{fig:infinPlatSig}
\end{figure}
\begin{figure}[htbp]
    \centering
    \subfloat[]{\label{fig:meEffKI}
    \includegraphics[width=2.5in]{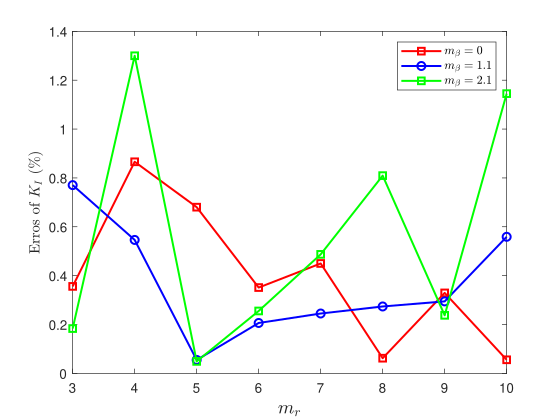}} \subfloat[]{\label{fig:meEffKII}
    \includegraphics[width=2.5in]{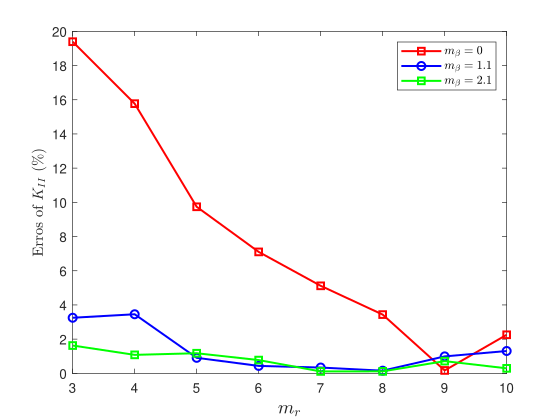}}
    \caption{Relative SIFs' error for $\theta=30^{\circ}$ at various integration contour and various $\beta$-PD element: (a) $K_I$, (b) $K_{II}$.}
    \label{fig:meEffect}
\end{figure}

\begin{figure}[htbp]
    \centering
    \includegraphics[width=3in]{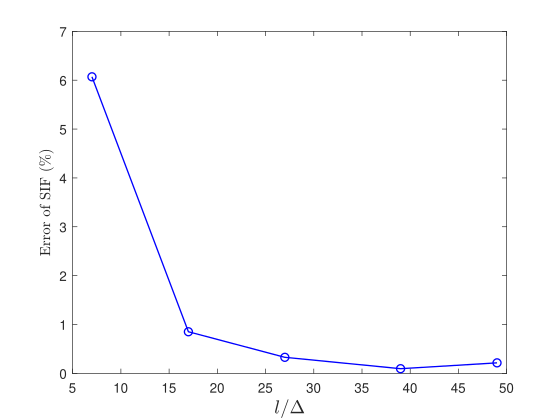}
    \caption{ Relative errors of SIF for model I loading with various mesh size.}
    \label{fig:meshRefine}
\end{figure}

\begin{figure}
    \centering
     \subfloat[]{\label{fig:diag_plat_gemo}
    \includegraphics[width=2.2in]{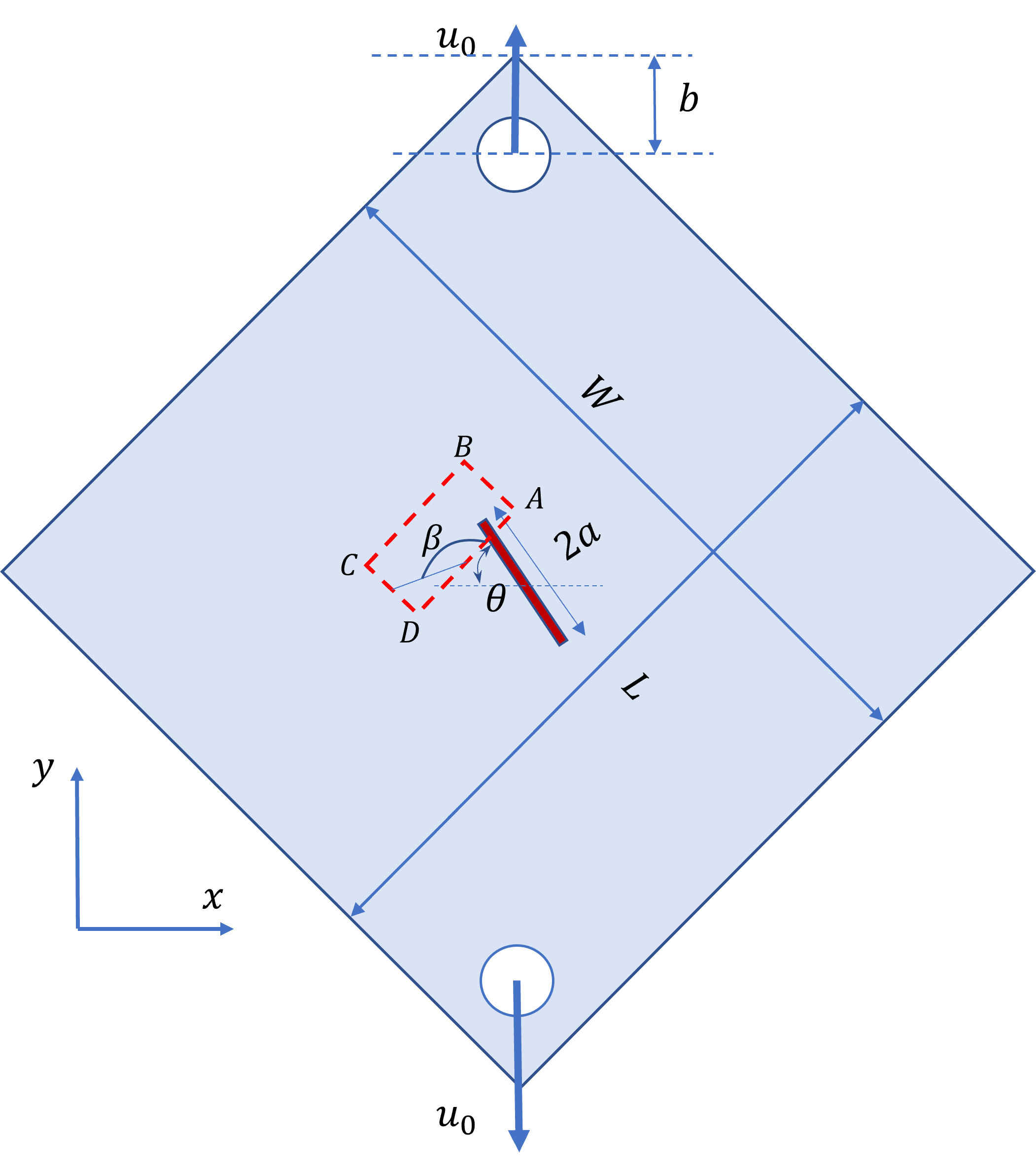}} 
    \subfloat[]{\label{fig:DiagMSH}
    \includegraphics[width=3.5in]{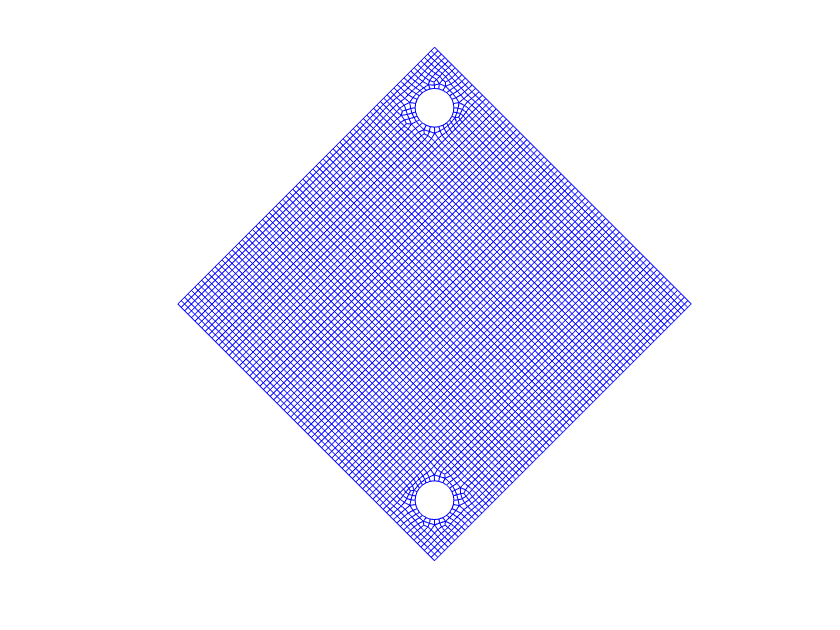}}
    \caption{A 2-D diagonal plate with an inclined crack under mixed-mode loading: (a) Geometry and boundary conditions, (b) Meshes.}
    \label{fig:diag_plat_gemo_mesh}
\end{figure}
\begin{figure}
    \centering
     \subfloat[]{\label{fig:DiagPlatUx}
    \includegraphics[width=2.5in]{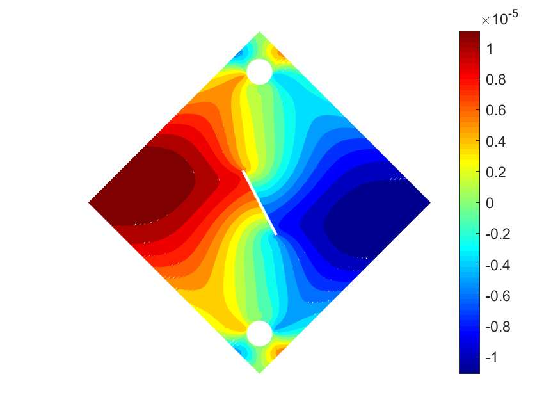}}  
    \subfloat[]{\label{fig:DiagPlatUy}
    \includegraphics[width=2.5in]{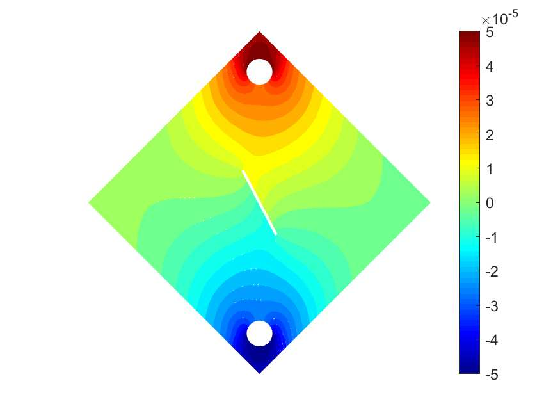}}
    \caption{Displacement results of the diagonal plate by adaptive PDLSM-FEM: (a) $u_x$, (b) $u_y$.}
    \label{fig:diag_plat_disp}
\end{figure}
\begin{figure}
    \centering
    \includegraphics[width=4in]{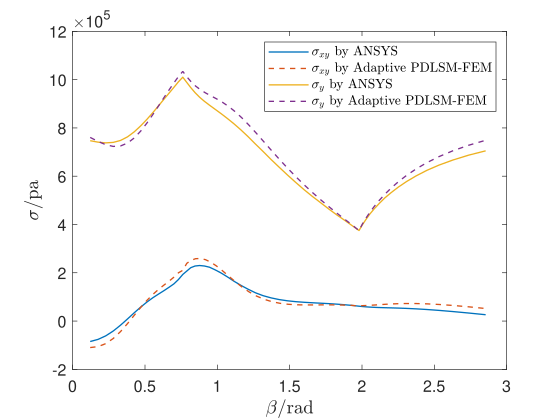}
    \caption{Comparison of FEM (ANSYS) against adaptive PDLSM-FEM for stress results along the path $A-B-C-D$ of the diagonal plate.}
    \label{fig:DiagPlatComS}
\end{figure}

 \begin{figure}[htbp]
    \centering
     \subfloat[]{\label{fig:DiagPlatSig_8}
    \includegraphics[width=2.5in]{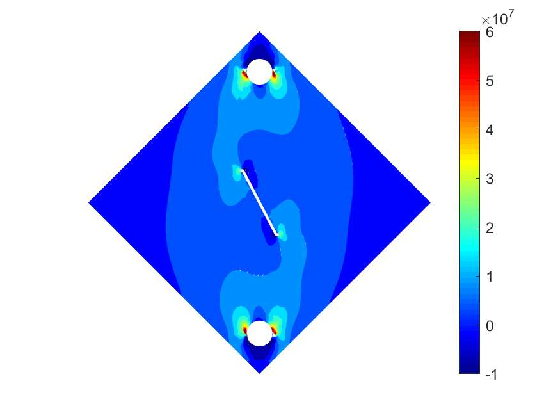}}  
    \subfloat[]{\label{fig:DiagPlatSig_22}
    \includegraphics[width=2.5in]{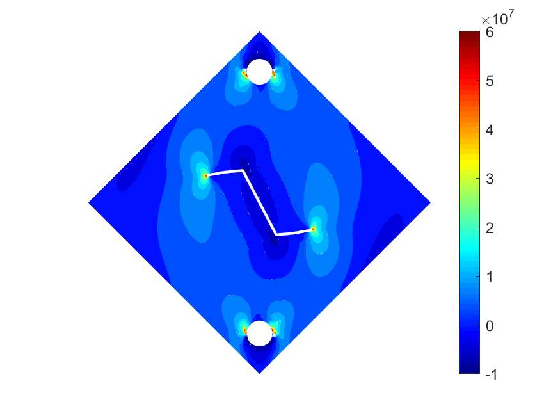}}\\
     \subfloat[]{\label{fig:DiagPlatSig_36}
    \includegraphics[width=2.5in]{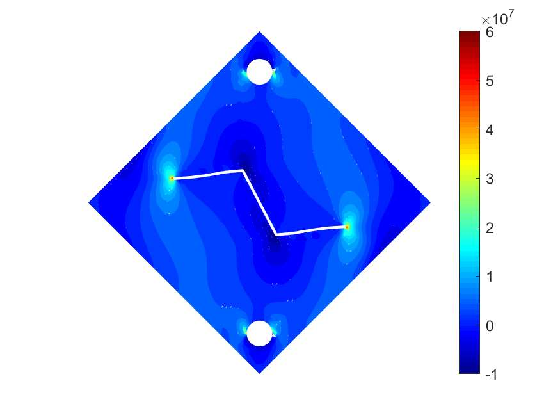}}
     \subfloat[]{\label{fig:DiagPlatSig_50}
    \includegraphics[width=2.5in]{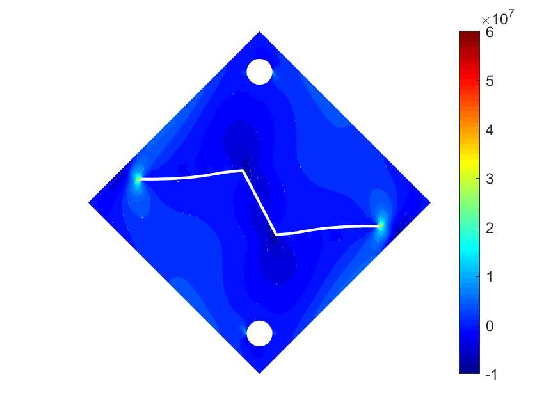}}
    \caption{Stress $\sigma_y$ contour of the 2-D diagonal plate during crack growth. }
    \label{fig:DiagS}
\end{figure}
\begin{figure}[htbp]
    \centering
    \includegraphics[width=3in]{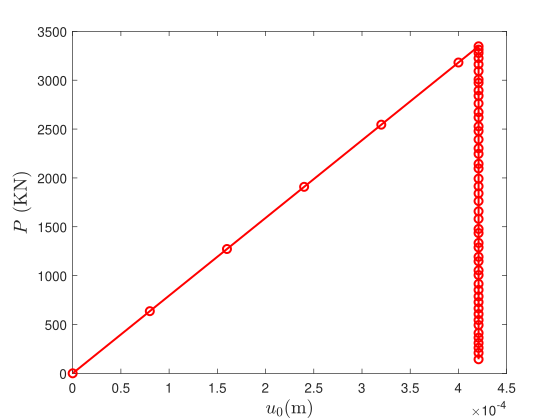}
    \caption{The reaction force of the 2-D diagonal plate during crack growth.}
    \label{fig:DiagPlatRF}
\end{figure}
\begin{figure}[htbp]
    \centering
    \includegraphics[width=4in]{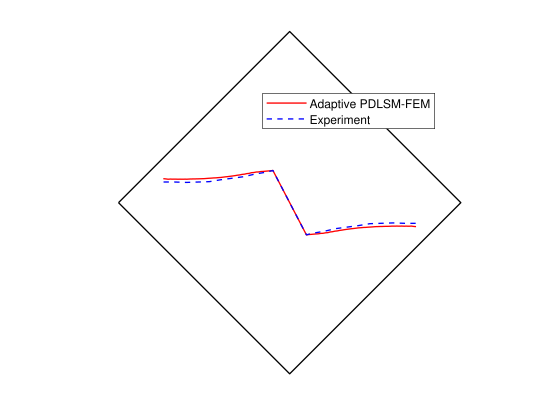}
    \caption{Crack propagation path of the 2-D diagonal plate by adaptive PDLSM-FEM and experiment \citep{ayatollahi_analysis_2009}.}
    \label{fig:DiagPath}
\end{figure}

\subsection{A 2-D infinite plate}

In this example, we study the SIFs of the 2-D infinite plate to validate the accuracy of the presented adaptive PDLSM-FEM. The effects of PD elements and the integral path of $I$-integral are studied. The crack remains stationary, and the primary purpose of this simulation is to investigate the model accuracy and effects of the $I$-integral contour and the $\beta$-PD element factor. 

\cref{fig:infiniPlatGemo} presents a 2-D square plate with a small crack subjected to uniform far-field tension $\sigma=70$ MPa. The plate's material properties are $E=70$ GPa and $\nu=0.33$. The plate's dimension is $1\times 1\ \mathrm{m}^2$, and the initial crack length is specified as $2a=40$ mm. Because of the sufficiently large ratio $W/a=50$, this plate can be approximately treated as infinity. A small sub-domain of dimension $0.1 \times 0.1\ \mathrm{m}^2$ at the plate center is finely discretized into 1521 elements, and the remaining part is coarsely meshed into 1728 elements ( \cref{fig:infiniPlatMSH}). For this example, the analytical solution of SIFs is defined as
\begin{equation}
    \begin{cases}
     K_I=\sigma\sqrt{\pi a} \cos^2{\beta},\\
     K_{II}=\sigma\sqrt{\pi a} \cos{\beta}\sin{\beta}.
    \end{cases}
\end{equation}

In this example, the based circle for the $I$-integral contour is selected as $r=m_r \Delta_{min}$, where $m_r=6$ and $\Delta_{min}$ is the minimum element size, and the horizon is $\delta_{(i)} = 3 \Delta_{(i)}$. The $\beta$-PD elements are defined as $m_{\beta}=2.1$.

\cref{fig:SIFs} presents the results of SIFs, $K_I$ and $K_{II}$, at various crack inclined angles $\theta$ from theoretical solution and numerical simulation. It is evident that the agreement between numerical and theoretical results is remarkable, and the relative error can be almost neglected. \cref{fig:infinPlatSig} shows the stress distributions around the crack by adaptive PDLSM-FEM with the crack inclined $\theta=40^\circ$. The result confirms that stresses concentrate at the crack tips.

The effects of integration radius $r=m_r \Delta_{min}$ and the $\beta$-PD element factor $m_{\beta}$ on the evaluation of SIFs are also studied for the plate with the crack inclination angle of $\theta=30^\circ$. The results are presented in \cref{fig:meEffect}. \cref{fig:meEffKI} shows that the relative errors of $K_I$ remain consistently below 1.3\% for various $m_{\beta}$ and $m_r$. Thus, the effect of the $\beta$-PD element factor $m_{\beta}$ on the evaluation of $K_I$ is practically negligible.

\cref{fig:meEffKII} shows the results of $K_{II}$. For $m_r\leq8$, the relative errors of $K_{II}$ are quite large, with the significant error being up to 19.4\% when there are no $\beta$-PD elements ($m_{\beta}=0$). In the presence of $\beta$-PD elements, the relative errors of $K_{II}$ are less than 3.25\% for $m_{\beta}=1.1$, and less than 1.6\% for $m_{\beta}=2.1$ at different $m_r$. Overall, the presence of $\beta$-PD elements has a significant beneficial effect on the evaluation of $K_{II}$. Although the $m_{\beta}$ has a negligible effect on the evaluation of $K_I$ for this case, the authors believe its effect on evaluating $K_I$ may not be neglected for other cases. With $\beta$-PD elements, the relative errors of SIFs are almost negligibly minor at different $m_r$. Thus, the evaluation of SIFs by $I$-integral is path independent. Although there is oscillation when $m_r$ is increased, we believe it is convergent because the relative error is negligible (In \cref{fig:meEffKI}, the maximum error of $K_I$ is around 1.3\%). 
In \cref{fig:meEffKII}, for $m_{\beta}=1.1$ and $m_{\beta}=2.1$, the error of $K_{II}$ does not change too much; however, it is evident that for all these $m_r$ values, the error of $K_{II}$ is minimal already, and it reaches the minimum error of $K_{II}$ when $m_r$ equals 8. It needs to be noted that when $m_{\beta}=0$, the error of $K_{II}$ keeps dropping as $m_r$ increases, and it reaches the minimum value when $m_r$ equals 9.0.
The effect of mesh refinement is investigated for this example, and the results are shown in \cref{fig:meshRefine}, in which $l=0.1$ m is the square's length with fine meshes (see \cref{fig:infiniPlatMSH}), and $\Delta$ is the mesh size around the crack. \cref{fig:meshRefine} shows the convergence with decreasing mesh size.

\subsection{A 2-D diagonal plate}

\begin{figure}
    \centering
     \subfloat[]{\label{fig:CTGeo}
    \includegraphics[width=2.2in]{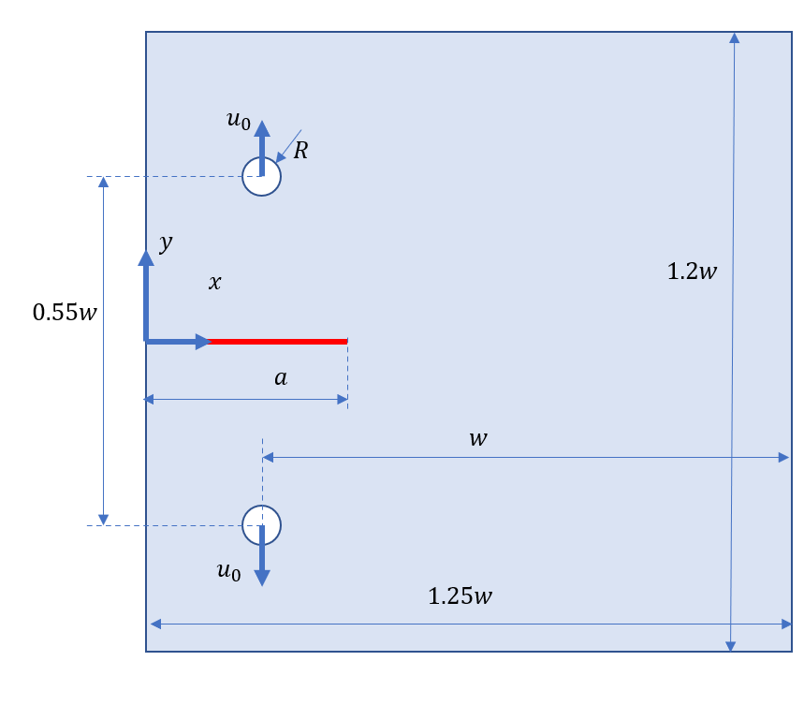}}
    \subfloat[]{\label{fig:CTmsh}
    \includegraphics[width=2.7in]{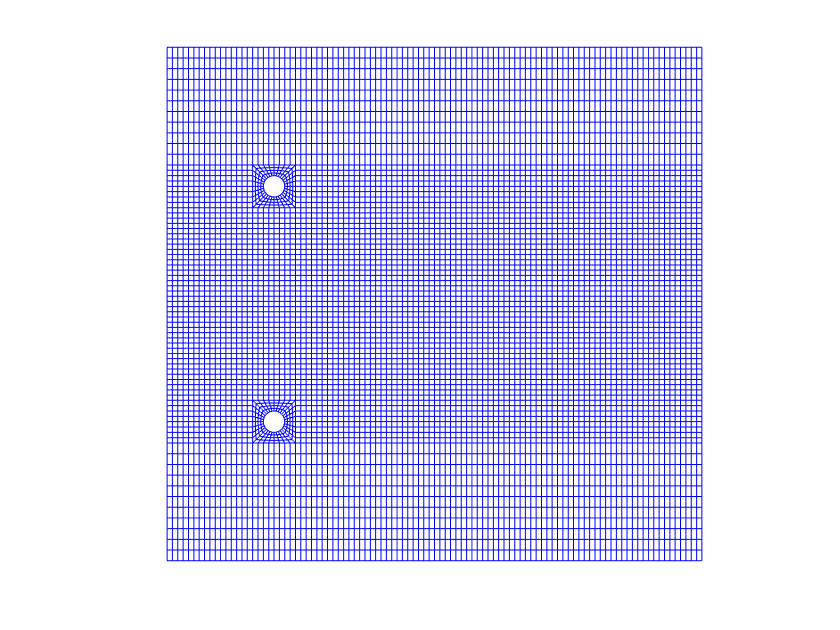}}
    \caption{Compact tension test: (a) Geometry and boundary conditions, (b) meshes.}
    \label{fig:CT_gemo_mesh}
\end{figure}

 \begin{figure}
    \centering
     \subfloat[]{\label{fig:CTSig4}
    \includegraphics[width=2.5in]{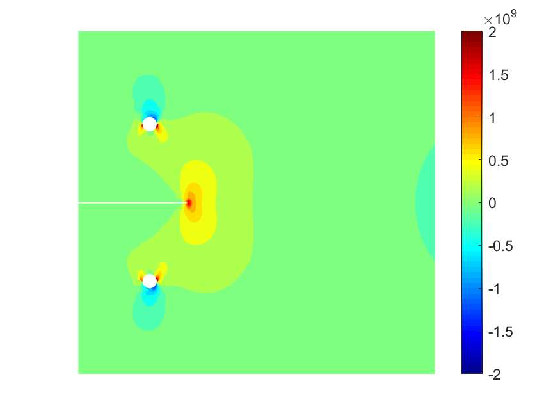}}  
    \subfloat[]{\label{fig:CTSig56}
    \includegraphics[width=2.5in]{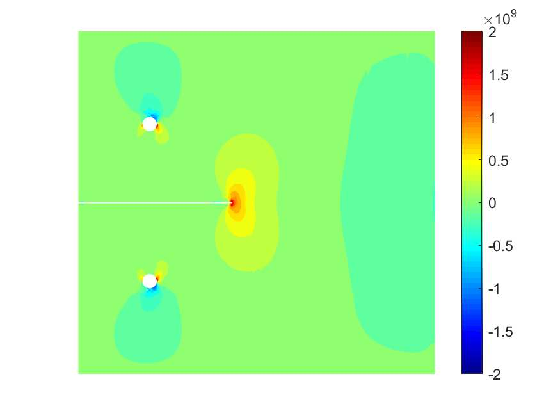}}\\
     \subfloat[]{\label{fig:CTSig108}
    \includegraphics[width=2.5in]{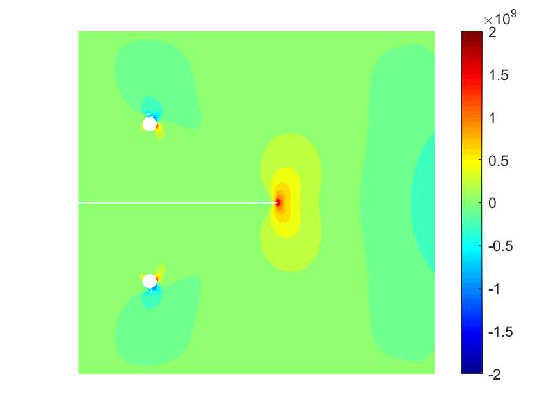}}
     \subfloat[]{\label{fig:CTSig160}
    \includegraphics[width=2.5in]{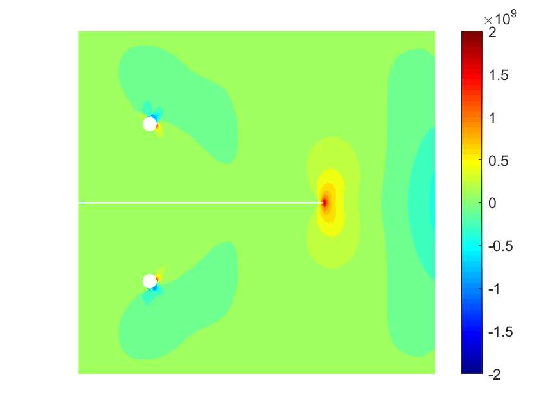}}
    \caption{Stress distributions of $\sigma_y$ of compact tension test: (a) $u_0=0.0476$ mm; (b) $u_0=0.0670$ mm; (c) $u_0=0.0928$ mm; (d) $u_0=0.1218$ mm.}
    \label{fig:CTS}
\end{figure}

\begin{figure}
    \centering
    \includegraphics[width=3in]{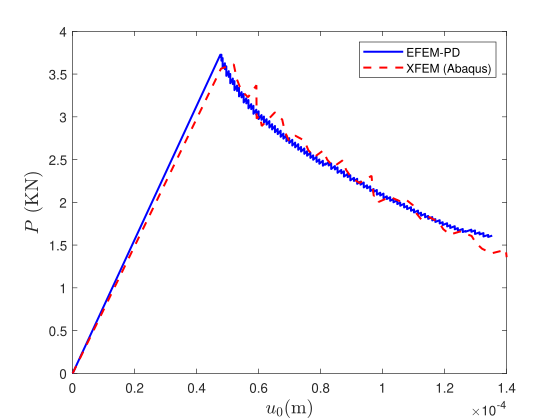}
    \caption{The reaction force the CT specimen under displacement control.}
    \label{fig:CTRF}
\end{figure}

For the second example, we first conduct a static simulation of the cracked body by both adaptive PDLSM-FEM and FEM (ANSYS) and compare the results to show the accuracy of the proposed model. We then simulate the quasi-static crack propagation using adaptive PDLSM-FEM and compare the simulation results with published experiments.

The problem to be simulated is a diagonal square plate with a pre-existing crack. As presented in \cref{fig:diag_plat_gemo}, the plate with dimension $150 \times 150 \times 5\ \mathrm{mm}^3$ is under displacement loading $u_0$. The crack is inclined with angle $\theta=62.5^\circ$, and its length is $2a=45$ mm. The loading hole has a radius of 8 mm, and the length between the corner of the plate and the center of the hole is $b=25$ mm. Young’s modulus is $E=2.94$ GPa, Poisson ratio is $\nu$= 0.38, and fracture toughness is $K_{Ic}=1.33$ MPa. As presented in \cref{fig:DiagMSH}, the problem domain is meshed into 5607 elements and 5785 nodes. The based circle for the $I$-integral contour is $r=6\Delta_{min}$. The horizon is specified as $\delta_{(i)}=3\Delta_{(i)}$. The $\beta$-PD element factor is set to be $m_{\beta}=3$.

A static stimulation for this problem is first performed with displacement loading $u_0=5\times 10^{-5}$ m. \cref{fig:diag_plat_disp} depicts the displacement field of the whole domain obtained by the adaptive PDLSM-FEM. The figure reveals apparent displacement discontinuities across the crack. \cref{fig:DiagPlatComS} shows the static results of the shear stresses, $\sigma_{xy}$, and the normal stresses, $\sigma_y$, by adaptive PDLSM-FEM and FEM (ANSYS) along path $A-B-C-D$ defined in \cref{fig:diag_plat_gemo}, where the key points' coordinates are $A\ (7.43,28.6)$ mm, $B\ (-13.8,49.9)$ mm, $C\ (-49.9,13.8)$ mm and $D\ (-13.8,49.9)$ mm. As shown in \cref{fig:DiagPlatComS}, the results of these two methods are in good agreement.

After the static simulation comparison, a quasi-static analysis is performed to simulate the crack propagation. The amount of crack growth for each step is set to be $d_c=\Delta_{min}$. Distributions of stress $\sigma_y$ of the diagonal plate at various stages of crack growth are shown in \cref{fig:DiagS}. The figure reveals that the stresses concentrate at the crack tip. Variations of the reaction force at the loading points with crack propagation are shown in \cref{fig:DiagPlatRF}. The figure shows that the crack starts to grow unstably when displacement loading reaches $u_0=0.421$ mm with the reaction force decreasing as the crack grows. The critical reaction force at instability is $P=3347$ N from the adaptive PDLSM-FEM simulation, whose relative error is 4.85\%, compared to the experimental result of $P_{cr}=3192$ N \citep{ayatollahi_analysis_2009}. \cref{fig:DiagPath} presents the crack propagation path of the plate from adaptive PDLSM-FEM simulation and the experiment, and they are in good agreement.

\subsection{Compact tension test}

\begin{figure}[htbp]
    \centering
         \subfloat[]{\label{fig:blocCracGemo}
    \includegraphics[width=2.8in]{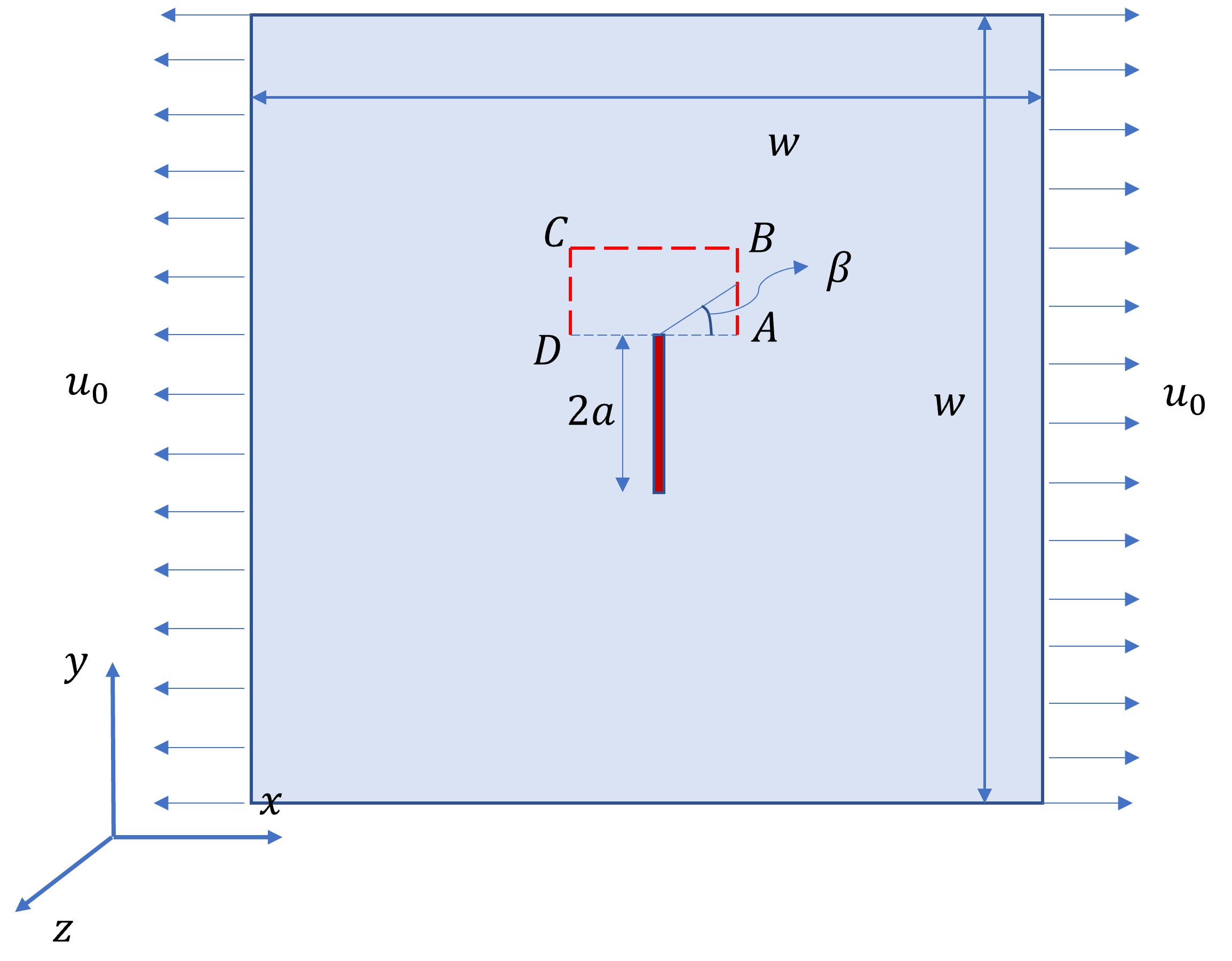}} \subfloat[]{\label{fig:blocCracMsh}
    \includegraphics[width=2.1in]{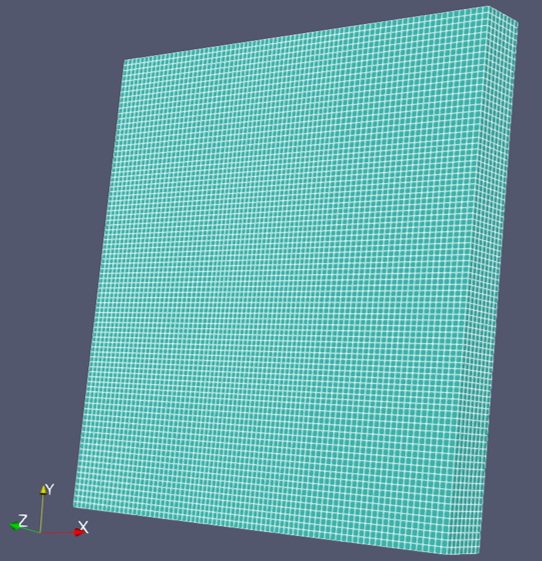}}
    \caption{The 3-D block with a pre-existing crack under displacement constraint: (a) Geometry and boundary conditions, (b) Hexahedron meshes.}
    \label{fig:blocCrac}
\end{figure}

\begin{figure}[htbp]
    \centering
         \subfloat[]{\label{fig:blocCracUx}
    \includegraphics[width=1.5 in]{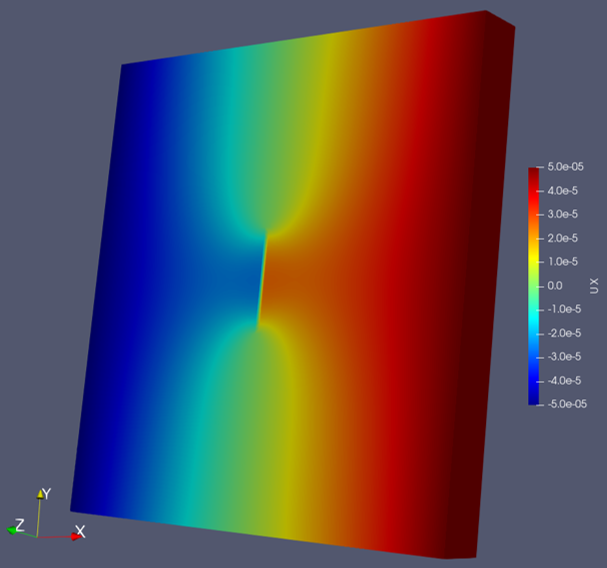}}\ \subfloat[]{\label{fig:blocCracUy}
    \includegraphics[width=1.5 in]{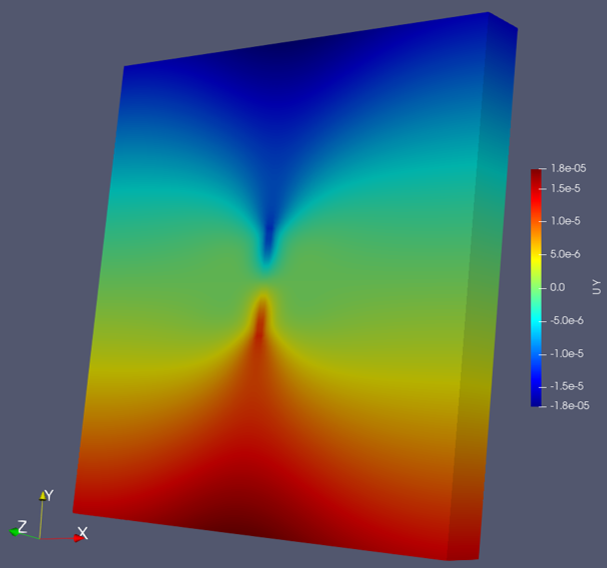}}\
    \subfloat[]{\label{fig:blocCracUz}
    \includegraphics[width=1.5 in]{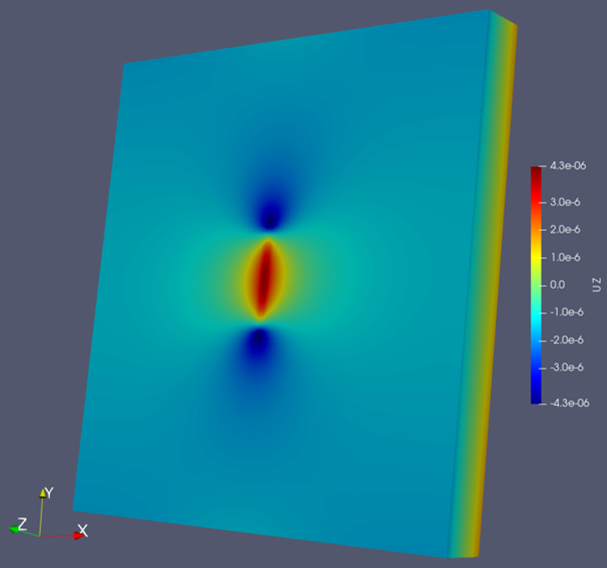}}
    \caption{Displacement solutions of the 3-D block with a pre-existing crack by adaptive PDLSM-FEM: (a) $u_x$, (b) $u_y$, (b) $u_z$.}
    \label{fig:blocCracU}
\end{figure}

\begin{figure}[htbp]
    \centering
     \subfloat[]{\label{fig:blocCracSx}
    \includegraphics[width=1.5 in]{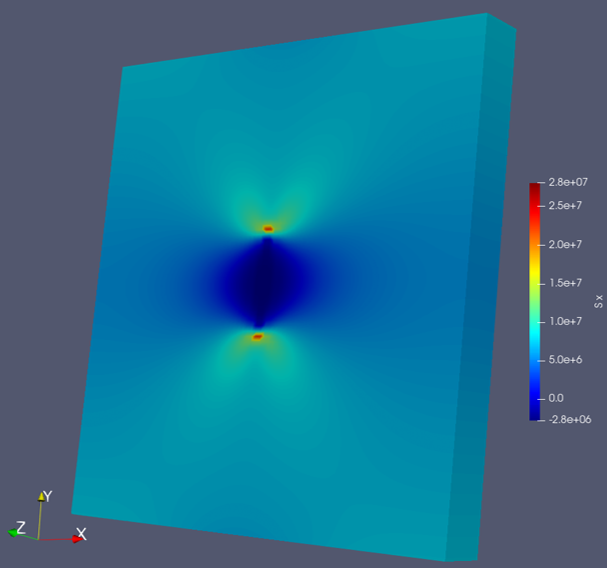}}\ \subfloat[]{\label{fig:blocCracSy}
    \includegraphics[width=1.5 in]{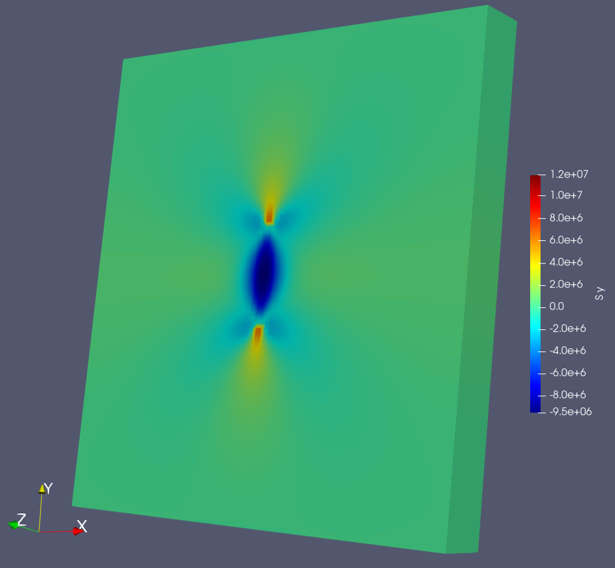}}\
    \subfloat[]{\label{fig:blocCracSxy}
    \includegraphics[width=1.5 in]{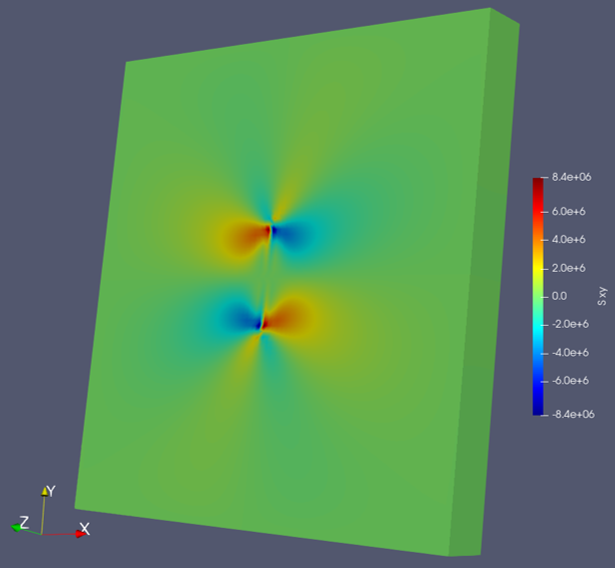}}
    \caption{Stress solutions of the 3-D block a with pre-existing crack by adaptive PDLSM-FEM: (a) $\sigma_x$, (b) $\sigma_y$, (c) $\sigma_{xy}$}
    \label{fig:blocCracS}
\end{figure}

\begin{figure}[htbp]
    \centering
         \subfloat[]{\label{fig:bloCracCompU}
    \includegraphics[width=2.3in]{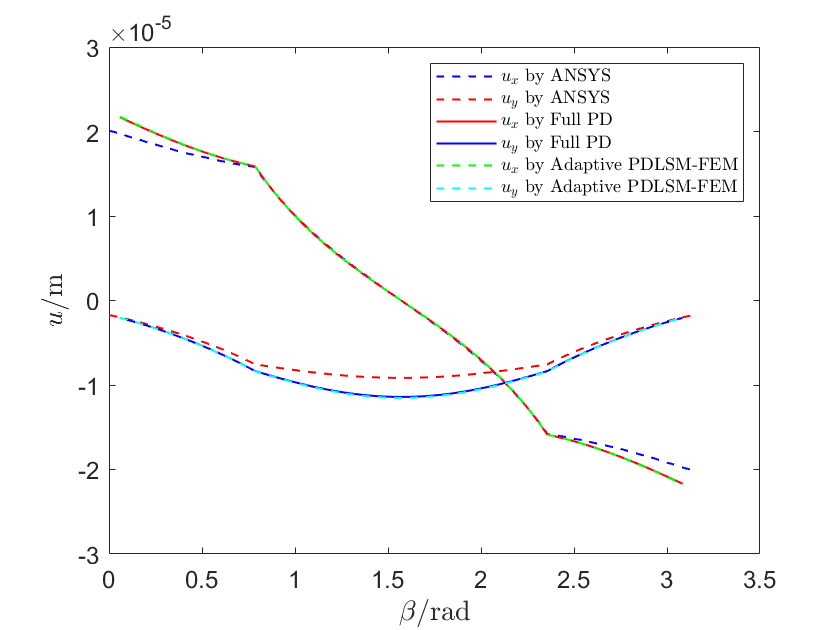}}\quad  \subfloat[]{\label{fig:bloCracCompS}
    \includegraphics[width=2.3 in]{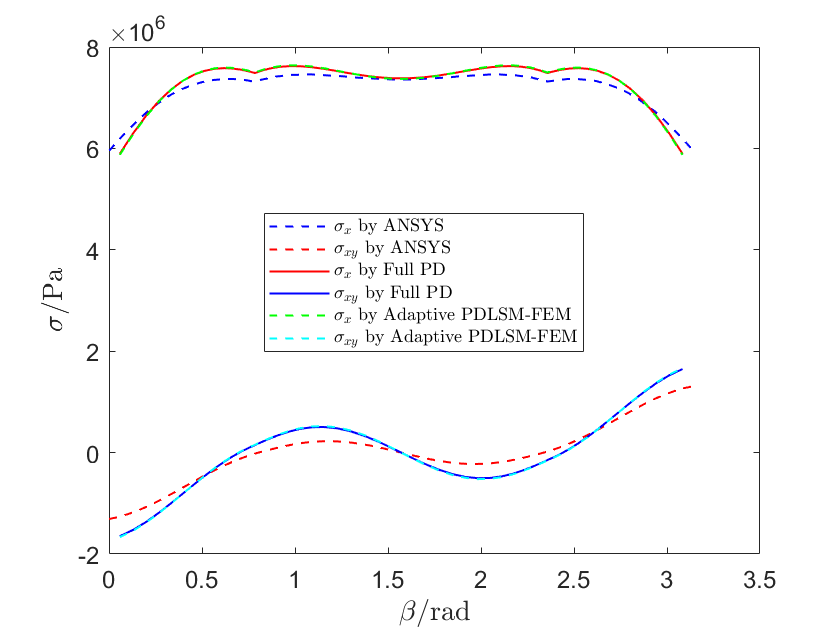}}
    \caption{Comparison of FEM (ANSYS) against adaptive PDLSM-FEM for displacement and stress results along with the path $A-B-C-D$ of the 3-D block: (a) Displacements Comparison; (b) Stresses comparison.}
    \label{fig:blocCracComp}
\end{figure}

In this example, the quasi-static crack propagation in a compact tension (CT) specimen is simulated. The results are compared with those obtained using the method of XFEM from ABAQUS. 

The compact tension (CT) specimen is subjected to displacement loading $u_0$, as depicted in \cref{fig:CTGeo}. The plate's dimension is specified as $w=50$ mm with thickness $B=1$ mm. The initial length of the crack is $a=15$ mm. The material properties of the plate are specified as $E=214 \ \mathrm{GPa}$, $\nu=0.27$, and $K_{Ic}=64.2\ \mathrm{MPa \sqrt{m}}$. As shown in \cref{fig:CTmsh}, the plate is meshed into 7834 nodes and 7628 elements. The amount of crack growth is defined as $d_c=\Delta_{min}$, where $\Delta_{min}$ is the minimum element size. The horizon is $\delta_{(i)}=3\Delta_{(i)}$, and the $I$-integral contour is defined as $r=6\Delta_{min}$. The $\beta$-PD elements are defined by $m_{\beta}=3$. 

\cref{fig:CTS} depicts the crack growth along with stress contours of $\sigma_y$ of the plate at various simulation steps. It reveals that the crack grows horizontally as expected and confirms the stresses concentrating at the crack tip with the tip stress fields, consistent with analytical solutions from LEFM. For comparison, the problem is also solved by XFEM with the virtual crack closure technique using ABAQUS. The variations of reaction forces with prescribed displacement $u_0$ of the CT test obtained from adaptive PDLSM-FEM and XFEM are plotted in \cref{fig:CTRF}. The figure shows a good agreement between the two methods.

\subsection{A 3-D block}

Because of the nonlocal nature of nodal interactions, PD simulations are typically much more computationally intensive than FEM simulations for similar problems, especially for 3-D cases. In this section, simulations of a 3-D body with a stationary crack using adaptive PDLSM-FEM, full PD in which PDLSM equations are applied to the entire elements, and FEM (using ANSYS) models are conducted and reported. The 3-D block is loaded under displacement control. The purpose of these simulations is to compare the model's accuracy and efficiency.

This 3-D example is much more computationally intensive than other previous 2-D examples. The adaptive PDLSM-FEM and full PD simulations are parallel run on four cores of the \emph{Moles} computing nodes of the high-performance cluster \emph{Beocat} at \emph{Kansas State University} by message passing interface. 

\cref{fig:blocCracGemo} shows the model geometry with a pre-existing crack subjected to displacement loading $u_0=5\times 10^{-5}$ m. The dimensions of the block are thickness $B=0.1$ m and width $w=1$ m. The crack length is $2a=0.2$ m. As shown in \cref{fig:blocCracMsh}, the whole domain is constructed by 52976 hexahedron elements and 61074 nodes. The horizon is specified as $\delta_{(i)}=3\Delta_{(i)}$. $\beta$-PD elements are defined by $m_{\beta}=2.1$.

\begin{table}[!htp]
\caption{\label{tab:perfEval}Wall time to assemble the global stiffness matrix and solve for displacement solutions, and the memory usages to store the global stiffness matrix, of the 3-D block with a pre-existing crack.}
\centering
\begin{tabular}{c|c|c|c|c}
\hline
~&\multicolumn{2}{c|}{Wall time (s)}&\multicolumn{2}{c}{Memory usage (MB)}\\\hline
$m_{\delta}$ & Adaptive PDLSM-FEM & Full PD &Adaptive PDLSM-FEM & Full PD\\\hline
3&	13.5&	63.7&	200.4&	700.6\\
4&	18.6&	169.3&	218.5&	1576.8\\
5&	25.6&	410.5&	247.1&	2940.9\\
\hline
\end{tabular}
\end{table}

\cref{fig:blocCracU} shows the $u_x$, $u_y$, and $u_z$ displacement fields, and \cref{fig:blocCracS} shows the $\sigma_x$, $\sigma_y$, and $\sigma_{xy}$ stress fields of the 3-D block from adaptive PDLSM-FEM. As revealed in these figures, stresses concentrate at the crack tips. \cref{fig:blocCracComp} compares the displacements $u_x$ and $u_y$ and the  stresses $\sigma_{xy}$ and $\sigma_x$, from adaptive PDLSM-FEM, full PD, and FEM (ANSYS), along path A-B-C-D as marked in  \cref{fig:blocCracGemo}, in which the key points' coordinates are $A\ (0.15,0.1,0)$ m, $B\ (0.15,0.25,0)$ m, $C\ (-0.15,0.25,0)$ m, and $D\ (-0.15,0.1,0)$ m. The comparison shows that the results from the full PD and from adaptive PDLSM-FEM are in excellent agreement and practically indistinguishable. It is observed that for the horizontal segment $B-C$, the results of $u_x$ from adaptive PDLSM-FEM and full PD model are close to those from ANSYS, while for the vertical segments $A-B$ and $C-D$, the $u_x$ from adaptive PDLSM-FEM and full PD model are slightly larger than those from ANSYS. Nevertheless, for the horizontal segment $B-C$, the $u_y$ from adaptive PDLSM-FEM and full PD model are slightly larger than those from ANSYS, while for the vertical segments $A-B$ and $C-D$, the $u_y$ from adaptive PDLSM-FEM and full PD model are very close to those from ANSYS. It is evident that for both stresses and displacements, the results by all three methods are in good agreement. 

To demonstrate the improved computational efficiency of the presented adaptive PDLSM-FEM, the wall time to assemble the global stiffness matrix and solve displacement solutions on four computing cores is tracked, and the memory usage for storing the global stiffness matrix for each core is monitored. As shown in \cref{tab:perfEval}, the wall time for adaptive PDLSM-FEM is 13.5-25.6 s for various horizon sizes $\delta_{(i)}=m_{\delta}\Delta_{(i)}$ ($m_{\delta}=3,4,5$), while the full PD model requires 63.7-410.5 s, which is 4.7-16.0 times that of the adaptive PDLSM-FEM model. The memory usage of adaptive PDLSM-FEM is 200.4-247.1 MB, while the full PD model requires 700.6-2940.9 MB, which is 3.5-11.9 times that of the adaptive PDLSM-FEM model. Therefore, the adaptive PDLSM-FEM dramatically improves the computational cost in CPU run time and memory usage compared to the full PD model.

\section{Conclusion}\label{sec:conc}

In this work, we have developed a framework of the adaptive PDLSM-FEM for modeling discontinuities analogous to the local enrichment of XFEM. In adaptive PDLSM-FEM, the Finite Element Method (FEM) is coupled with Peridynamics (PD) with adaptivity to minimize the PD region to maximize the computational efficiency. With the framework of adaptive PDLSM-FEM, the elements intersecting with the crack line and their neighboring elements are governed by the PDLSM model that is capable of modeling the progression of bond breakages and crack growth. Conventional FEM governs the remaining elements. The global stiffness matrix and governing equations of the whole problem domain in 2-D and 3-D are derived. A numerical procedure for the $I$-integral to calculate SIFs for 2-D problems is proposed and implemented. Crack propagation is modeled using quasi-static with maximum hoop tension stress criterion. New contributions of this work include adaptivity coupling PDLSM-FEM for minimizing the PD region and the application of the adaptive PDLSM-FEM to quasi-static crack propagation analysis.

Several numerical examples are conducted to verify the efficacy of the proposed model. The remarks of numerical simulations are as below:
\begin{enumerate}
    \item One 2-D infinite plate under mixed-mode loading is analyzed, and its SIFs are calculated. Comparing the SIFs between analytical and numerical solutions reveals excellent accuracy of the presented adaptive PDLSM-FEM.
    \item Studying the evaluation of SIFs reveals that $I$-integral is path independent, and $\beta$-PD elements have a significant beneficial effect on evaluations of SIFs.
    \item Static simulation of a diagonal plate is performed, and the static solutions align with ANSYS results. 
    \item Adaptive PDLSM-FEM also simulates the quasi-static crack propagation of the 2-D diagonal plate under displacement control. The crack path and peak loading are in close agreement with experiment observation.
    \item Simulations of a compact tension test by adaptive PDLSM-FEM and XFEM (ABAQUS) are performed, and the reaction forces from both methods are in good agreement.
    \item The static solutions of a 3-D block with pre-existing crack are obtained by adaptive PDLSM-FEM, full PD, and FEM (ANSYS), and the comparison reveals that the results from the three methods are in excellent agreement.
    \item The performance evaluation of the 3-D example shows that adaptive PDLSM-FEM saves up to 16.0 times in terms of computing time and up to 11.9 times in terms of memory usage against the pure PD model.

\end{enumerate}




\bibliographystyle{elsarticle-num-names}  
\bibliography{cas-refs}
\end{document}